\documentstyle{amsppt}
\voffset-10mm
\magnification1200
\pagewidth{130mm}
\pageheight{204mm}
\hfuzz=2.5pt\rightskip=0pt plus1pt
\binoppenalty=10000\relpenalty=10000\relax
\TagsOnRight
\loadbold
\nologo

\let\ge\geqslant
\let\wt\widetilde
\let\wh\widehat
\let\[\lfloor
\let\]\rfloor
\redefine\d{\roman d}

\define\Li{\operatorname{Li}}
\define\Tr{\operatorname{Tr}}
\redefine\pmod#1{\;(\operatorname{mod}#1)}

\define\sect#1#2#3{\parshape=3 4mm 117mm 9mm 112mm 9mm 112mm
\noindent\eightpoint{\bf#1}\enspace#2\dotfill\dotfill
\noindent\vskip-\normalbaselineskip\parshape=1 4mm 122mm {}
\hfill\eightpoint#3}
\newcount\superscount
\def\supers{\global\advance\superscount1\relax\number\superscount}
\topmatter
\title
Differential equations, \\
mirror maps and zeta values
\endtitle
\author
Gert Almkvist {\rm(Lund) and}
Wadim Zudilin\footnotemark"$^\ddag$"\ {\rm(Moscow)}
\endauthor
\comment
\address
\hbox to70mm{\vbox{\hsize=70mm%
\leftline{Matematikcentrum}
\leftline{Lunds Universitet}
\leftline{Matematik MNF, Box 118}
\leftline{SE-221\,00 Lund, SWEDEN}
}}
\endaddress
\email
\tt gert\@maths.lth.se
\endemail
\address
\hbox to70mm{\vbox{\hsize=70mm%
\leftline{Department of Mechanics and Mathematics}
\leftline{Moscow Lomonosov State University}
\leftline{Vorobiovy Gory, GSP-2}
\leftline{119992 Moscow, RUSSIA}
}}
\endaddress
\email
\tt wadim\@ips.ras.ru
\endemail
\endcomment
\date
\hbox to100mm{\vbox{\hsize=100mm%
\centerline{E-print \tt math.NT/0402386}
\smallskip
\centerline{16 December 2004}
}}
\enddate
\abstract
The aim of this work is an analytic investigation of differential
equations producing mirror maps as well as giving new examples
of mirror maps; one of these examples is
related to (rational approximations to) $\zeta(4)$.
We also indicate certain observations that might become a subject
of further research.
\endabstract
\endtopmatter
\rightheadtext{Differential equations, mirror maps and zeta values}
\leftheadtext{G.~Almkvist and W.~Zudilin}
\footnotetext"$^\ddag$"{The work of the second author
is supported by an Alexander von Humboldt
research fellowship and partially supported by grant
no.~03-01-00359 of the Russian Foundation for Basic Research.}
\document

The existence of this paper is due to the following observations.
With Ap\'ery's proof of the irrationality of~$\zeta(2)$
and~$\zeta(3)$ followed certain 2nd- and 3rd-order linear differential
equations (see~\cite{Be1}--\cite{Be3}, \cite{BP}).
Doing a similar
construction for $\zeta(4)$ resulted in a 5th-order differential
equation. This equation was similar to the linear differential
equations occuring in Calabi--Yau theory (except the order was 5
instead of~4). Computing the Lambert series of an analogue
of the Yukawa coupling we got integer coefficients divisible
by the square of the degree (not by the cube as in the Calabi--Yau
case). Then we managed to pull back the 5th-order differential
equation to one of order~4, which had all the properties
of a Calabi--Yau equation. This part is the main objective
of Sections~1--4. Quite recently we found two more examples
relating fourth order differential equations to simultaneous
rational approximations to di- and trilogarithms
(see the end of Section~4).

Collecting known cases of the Calabi--Yau equations,
14~of those were classical hypergeometric and 15~remaining
were discovered by V.\,V\. Batyrev, D.~van Straten et~al
in~\cite{BS1}, \cite{BS2}, \cite{Str}, we tried to extend
the list by application of the algorithm of creative
telescoping due to Gosper and Zeilberger. This seems to
be a pure machinery (cf.\ Section~5 below), in which we cannot
predict any end.
We also tried to figure out possible algebraic transformations
between them and our new example related to $\zeta(4)$. We did not succeed
in this, but doing that we found two other Calabi--Yau equations
emanating from quadratic transformations
of hypergeometric Calabi--Yau solutions of type~${}_4F_3$.
This is the subject of Section~6. Section~7 describes
another (rather rich!) possibilities for deriving fake 4th-order
linear differential equations. Summarizing methods of Sections~6 and~7
resulted in new classes of transformations which we present in Section~8.

The divisibility properties of the coefficients
of the Lambert series for the Yukawa coupling~$K(q)$
are equivalent to high-order Kummer congruences
between the coefficients of the power series for~$K(q)$.
The latter coefficients can be realized as the number
of fix points of iterates of a set map $T\:X\to X$,
where the number of orbits of order~$n$ is divisible
by~$n^2$. We figure out these very natural relations
in Section~9.

Our final table (mostly represented
in the joint contribution \cite{AESZ, Table~A})
contains more than 200 cases of
Calabi--Yau equations (see also van Enckevort's
electronic database~\cite{En}, which provides
more advanced knowledge, like instanton numbers up to~20,
in some cases also the number of
elliptic curves, monodromy matrices etc.).
In almost all cases, the mirror
map $z(q)$ divided by~$q$ seems to be a high power
of a series with integer coefficients (cf\. \cite{LY}).
For this we have no explanation at all,
but we collect this experimental
knowledge in another table \cite{AESZ, Table~B}.
The last table \cite{AESZ, Table~C} presents a brief systematic
guide to the main table from \cite{AESZ, Table~A}, with entries
ordered according to instanton numbers%
\footnote{With the kind permit of our coauthors
D.~van Straten and C.~van Enckevort, we include
the last table from \cite{AESZ} in the version
of the present paper appearing in
the BIRS workshop volume
on Calabi--Yau varieties and mirror symmetry.
The E-print version of this article at
{\tt http://arXiv.org/\allowlinebreak math.NT/\allowlinebreak 0402386}
also contains Section~10 and proofs of Propositions~1 and~3.}.

Finally, in Section~10, we find a 6th-order linear
differential equation coming from the simultaneous
approximations to $\zeta(3)$ and $\zeta(5)$,
and also a similar 6th-order equation for the series
$\sum_{n=0}^\infty z^n\sum_{k=0}^n\binom nk^6$.
Both cases seem to admit arithmetic properties very
close to those for Calabi--Yau cases, but for the latter
series we have a new phenomenon---a free parameter.
This fact does not look apparent: we briefly indicate
some similar cases as well.

\subsubhead
Acknowledgements
\endsubsubhead
First of all we want to thank Duco van Straten, who
has shown great interest in our work and suggested the
use of Hadamard products. He and Christian van Enckevort
contributed with finding several new Calabi--Yau equations, usually
coming from geometry. Combining our examples and methods
with those from~\cite{ES} resulted in the joint project~\cite{AESZ}.
Frits Beukers gave us details of computations
on \cite{Be4} and inspiration to start this work; his comments
helped us a lot to improve the first version of our manuscript.
Helena Verrill communicated us a couple of new equations and
Arne Meurman showed an algorithm to find the Hadamard product.
Special gratitude is due to Noriko Yui who kindly proposed us
to contribute the BIRS workshop volume
on Calabi--Yau varieties and mirror symmetry.

\head
Contents
\endhead
\eightpoint

\sect{\phantom01.}{Maximal unipotent monodromy}{3}

\sect{\phantom02.}{$4$th-order linear differential equations}{3}

\sect{\phantom03.}{Wronskian formalism}{6}

\sect{\phantom04.}{Strange integrality related to $\zeta(4)$}{9}

\sect{\phantom05.}{Strategy of finding new Calabi--Yau equations}{12}

\sect{\phantom06.}{Quadratic transformations and mirror maps}{14}

\sect{\phantom07.}{Hadamard products}{18}

\sect{\phantom08.}{More transformations}{26}

\sect{\phantom09.}{Supercongruences and $k$-realizable series}{31}

\sect{10.}{$6$th- and higher order differential equations}{33}

\sect{\phantom{00.}}{References}{36}

\newpage

\tenpoint

\head
1. Maximal unipotent monodromy
\endhead

Consider a linear differential equation
$$
y^{(s)}+a_{s-1}(z)y^{(s-1)}+\dots+a_1(z)y'+a_0(z)y=0,
\tag1.1
$$
where the prime stands for the $z$-derivative and
$a_0,a_1,\dots,a_{s-1}$ are meromorphic functions in variable~$z$.
Assuming that $z=0$ is a regular singularity of
the differential equation~\thetag{1.1},
we may write its coefficients in the form
$$
a_{s-j}(z)=z^{-j}\tilde a_{s-j}(z),
\qquad j=1,\dots,s,
$$
where the functions $\tilde a_{s-j}(z)$ are analytic at $z=0$.
The roots of the {\it indicial equation\/}
$$
\lambda(\lambda-1)\dotsb(\lambda-s+1)
+\tilde a_{s-1}(0)\lambda(\lambda-1)\dotsb(\lambda-s+2)+\dots
+\tilde a_1(0)\lambda+\tilde a_0(0)
=0
$$
determine the {\it exponents\/} of the differential equation~\thetag{1.1}
at the point $z=0$. Following \cite{Mo}, we will say that
our differential equation~\thetag{1.1}
is of {\it maximal unipotent monodromy\/}
(briefly, MUM) if its exponents at $z=0$ are all zero.
The Frobenius method (see, e.g., \cite{In, Section~16.1})
gives one a constructive way for writing
a basis of solutions to \thetag{1.1} in the form
$$
\gathered
y(z,\rho)=\sum_{n=0}^\infty A(n,\rho)z^{n+\rho}\pmod{\rho^s}
=y_0(z)+y_1(z)\rho+\dots+y_{s-1}(z)\rho^{s-1},
\\
A(0,\rho)=1\pmod{\rho^s}
\endgathered
\tag1.2
$$
(in particular, $A(0,0)=y_0(0)=1$),
where by definition we set
$$
z^\rho\pmod{\rho^j}=1+\log z\cdot\rho+\frac{\log^2z}2\cdot\rho^2
+\dots+\frac{\log^{j-1}z}{(j-1)!}\cdot\rho^{j-1}.
$$
The components $y_0,y_1,\dots,y_{s-1}$ in~\thetag{1.2} define
the {\it Frobenius basis\/} of the MUM differential
equation~\thetag{1.1} at $z=0$.


\head
2. $4$th-order linear differential equations
\endhead

Consider a linear homogeneous differential equation of order~$4$,
$$
y^{(4)}+a_3(z)y^{(3)}+a_2(z)y''+a_1(z)y'+a_0(z)y=0,
\tag2.1
$$
where $a_0,a_1,a_2,a_3$ are rational functions in the variable~$z$.
Suppose that \thetag{2.1} is MUM and let $y_0,y_1,y_2,y_3$
be the Frobenius basis of solutions to~\thetag{2.1}.
Set $t=y_1/y_0$.

\proclaim{Proposition 1 \cite{Al2}}
In the above notation, the condition
$$
a_1=\frac12a_2a_3-\frac18a_3^3+a_2'-\frac34a_3a_3'-\frac12a_3''
\tag2.2
$$
is equivalent to any of the following two possibilities:
$$
\frac{\d^2}{\d t^2}(y_2/y_0)
=\frac{\exp\bigl(-\frac12\int a_3(z)\,\d z\bigr)}
{y_0^2\bigl(\frac{\d t}{\d z}\bigr)^3}
\quad\text{and}\quad
\frac{\d^2}{\d t^2}(y_3/y_0)
=t\frac{\d^2}{\d t^2}(y_2/y_0).
\tag2.3
$$
\endproclaim

\demo{Proof}
It can be easily checked that the functions
$$
T_1=\frac{y_1}{y_0}\;(=t), \quad
T_2=\frac{y_2}{y_0}, \quad
T_3=\frac{y_3}{y_0}
$$
satisfy the 4th-order linear differential equation
$$
T^{(4)}+b_3T^{(3)}+b_2T''+b_1T'=0,
\tag2.4
$$
where
$$
\gather
b_1=a_1+2a_2\frac{y_0'}{y_0}+3a_3\frac{y_0''}{y_0}+4\frac{y_0'''}{y_0},
\\
b_2=a_2+3a_3\frac{y_0'}{y_0}+6\frac{y_0''}{y_0},
\qquad
b_3=a_3+4\frac{y_0'}{y_0}.
\endgather
$$
Since
$$
\frac{\d}{\d t}=\biggl(\frac{\d t}{\d z}\biggr)^{-1}\frac{\d}{\d z}
=\frac1{T_1'}\,\frac{\d}{\d z},
$$
we have
$$
\frac{\d^2}{\d t^2}\biggl(\frac{y_2}{y_0}\biggr)
=\frac{\d^2T_2}{\d t^2}
=\frac{\d}{\d t}\biggl(\frac{T_2'}{T_1'}\biggr)
=\frac1{T_1'}\biggl(\frac{T_2''}{T_1'}
-\frac{T_2'T_1''}{(T_1')^2}\biggr)
=\frac R{(\d t/\d z)^3},
\tag2.5
$$
where
$$
R=T_1'T_2''-T_1''T_2'
=\vmatrix T_1' & T_2' \\ T_1'' & T_2'' \endvmatrix.
\tag2.6
$$
With the help of \thetag{2.4} we deduce
$$
\gathered
R'=\vmatrix T_1' & T_2' \\ T_1''' & T_2''' \endvmatrix,
\quad
R''=\vmatrix T_1'' & T_2'' \\ T_1''' & T_2''' \endvmatrix
-b_3R'-b_2R,
\\
R'''=-b_3(R''+b_3R'+b_2R)+b_1R-(b_3R')'-(b_2R)'.
\endgathered
\tag2.7
$$
Therefore, the function $R$ satisfies the 3rd-order differential equation
$$
R'''+c_2R''+c_1R'+c_0R=0,
\tag2.8
$$
where
$$
c_0=b_2b_3-b_1+b_2', \quad c_1=b_2+b_3^2+b_3', \quad c_2=2b_3.
$$
We now claim that the function
$$
\tilde R
=\frac1{y_0^2}\exp\Bigl(-\frac12\int a_3(x)\,\d x\Bigr)
$$
also satisfies the linear differential equation~\thetag{2.8}
if and only if \thetag{2.2} holds.
By differentiating we obtain
$$
\align
\tilde R'
&=-\biggl(\frac{a_3}2+2\frac{y_0'}{y_0}\biggr)\tilde R,
\\
\tilde R''
&=\biggl(\frac{a_3^2}4-\frac{a_3'}2+2a_3\frac{y_0'}{y_0}
-2\frac{y_0''}{y_0}+6\frac{(y_0')^2}{y_0^2}\biggr)\tilde R,
\\
\tilde R'''
&=\biggl(\frac{3a_3a_3'}4-\frac{a_3''}2-\frac{a_3^3}8
+3\Bigl(a_3'-\frac12a_3^2\Bigr)\frac{y_0'}{y_0}
+3a_3\frac{y_0''}{y_0}-2\frac{y_0'''}{y_0}
\\ &\qquad
-9a_3\frac{(y_0')^2}{y_0^2}
+18\frac{y_0'y_0''}{y_0^2}-24\frac{(y_0')^3}{y_0^3}\biggr)\tilde R,
\endalign
$$
hence a necessary and sufficient condition for the function $\tilde R$
to satisfy \thetag{2.8} is the required one:
$$
\biggl(\frac{a_2a_3}2-\frac{a_3^3}8+a_2'
-\frac{3a_3a_3'}4-a_1\biggr)\tilde R=0.
$$
It remains to compare the first three terms (i.e., coefficients of
$x^{-3}$, $x^{-2}$, and $x^{-1}$)
in order to verify the equality $\tilde R=R$.

We are now required to check the second equality in~\thetag{2.3}
provided \thetag{2.2} holds.
Changing $R$ to $R_1$ and $T_2$ to $T_3$ in~\thetag{2.6}, \thetag{2.7}
we obtain that the function $R_1$ satisfies the differential
equation~\thetag{2.8}. On the other hand, differentiating the product
$T_1R=tR$ with the help of the differential equation~\thetag{2.4}
for $T_1$ and
of~\thetag{2.8} for $R$, we deduce that the function $tR$ also
satisfies \thetag{2.8}. Comparing of the first three coefficients
again yields $R_1=tR$, hence formulas \thetag{2.5} and
$$
\frac{\d^2}{\d t^2}\biggl(\frac{y_3}{y_0}\biggr)
=\frac{\d^2T_3}{\d t^2}
=\frac{\d}{\d t}\biggl(\frac{T_3'}{T_1'}\biggr)
=\frac{R_1}{(\d t/\d z)^3}
$$
complete the proof.
\qed
\enddemo

If the differential equation \thetag{2.1} is MUM and
assumption \thetag{2.2} holds, we define the {\it mirror map\/}
$z=z(q)$ as the inverse of $q(z)=e^t=\exp(y_1/y_0)$
and the {\it Yukawa coupling\/} as
$$
K(q)=N_0\cdot\frac{\d^2}{\d t^2}(y_2/y_0),
$$
where $N_0$ is some integer different from~$0$.
The interesting ({\it Calabi--Yau\/}) cases
are those when $z(q)\in\Bbb Z[[q]]$
and
$$
K(q)=\sum_{n=0}^\infty C_nq^n
=N_0+\sum_{l=1}^\infty\frac{N_ll^3q^l}{1-q^l}
\tag2.9
$$
with
$$
N_0=C_0\in\Bbb Z \qquad\text{and}\qquad
N_l=\frac1{l^3}\sum_{d\mid l}\mu\Bigl(\frac ld\Bigr)C_d
\in\Bbb Z
\quad\text{for $l=1,2,\dots$}
$$
and $\mu$ is the M\"obius function
(and, usually, $N_l>0$).

\proclaim{Proposition 2}
Under hypothesis \thetag{2.2} we have the identity
for Wronskian determinants:
$$
\vmatrix y_0 & y_3 \\ y_0' & y_3' \endvmatrix
=\vmatrix y_1 & y_2 \\ y_1' & y_2' \endvmatrix.
$$
\endproclaim

\demo{Proof}
The required identity is equivalent to
$y_0^2(y_3/y_0)'=y_1^2(y_2/y_1)'$, i.e.
$$
\biggl(\frac{y_3}{y_0}\biggr)'
=\biggl(\frac{y_1}{y_0}\biggr)^2
\biggl(\frac{y_0}{y_1}\cdot\frac{y_2}{y_0}\biggr)'
=t^2\biggl(\frac1t\,\frac{y_2}{y_0}\biggr)'.
$$
Since $\frac{\d}{\d z}=\frac{\d t}{\d z}\cdot\frac{\d}{\d t}$,
verification of the last identity is equivalent to showing that $Y=0$,
where
$$
Y=\frac{\d}{\d t}\biggl(\frac{y_3}{y_0}\biggr)
-t^2\frac{\d}{\d t}\biggl(\frac1t\,\frac{y_2}{y_0}\biggr)
=\frac{\d}{\d t}\biggl(\frac{y_3}{y_0}\biggr)
+\frac{y_2}{y_0}-t\frac{\d}{\d t}\biggl(\frac{y_2}{y_0}\biggr).
$$
But
$$
\frac{\d Y}{\d t}
=\frac{\d^2}{\d t^2}\biggl(\frac{y_3}{y_0}\biggr)
-t\frac{\d^2}{\d t^2}\biggl(\frac{y_2}{y_0}\biggr)
=0
$$
by Proposition~1. Therefore $Y=\operatorname{const}$, and
verifying near $z=0$ we get $Y=0$ (all $\log$-terms disappear).
This completes the proof.
\qed
\enddemo

\head
3. Wronskian formalism
\endhead

Consider an arbitrary pair $y(z),\tilde y(z)$ of linearly independent
solutions of a 4th-order differential equation~\thetag{2.1}
(not necessarily MUM).
Define the Wronskian determinant
$$
Y=\vmatrix y & \tilde y \\ y' & \tilde y' \endvmatrix.
$$

\proclaim{Proposition 3}
The function $Y=Y(z)$ satisfies a $6$th-order linear differential equation
with coefficients depending only on the coefficients of~\thetag{2.1}
\rom(i.e., independent of the choice of the pair $y,\tilde y$\rom),
called the exterior square of~\thetag{2.1}.
Moreover, this $6$th-order differential equation
reduces to a $5$th-order differential equation
\rom(i.e., the coefficient of $Y^{(6)}$ is zero\rom) if and only if
condition \thetag{2.2} holds.
\endproclaim

\remark{Remark \rom1}
In other words, condition \thetag{2.2} is equivalent
to the exterior square having order~5. As pointed out by Beukers,
by this equivalence properties~\thetag{2.3} are the derived
properties which are handy to calculate the Yukawa coupling.
Nevertheless, we kept the original proof of Proposition~1
to make it independent of considerations in this section.
\endremark

\remark{Remark \rom2}
The clain of Proposition~3 may be proved in {\tt Maple}
by the {\tt DEtools} command {\tt exterior\_power}.
\endremark

\demo{Proof}
We have
$$
\gather
Y'=\vmatrix y & \tilde y \\ y'' & \tilde y'' \endvmatrix,
\qquad
u_3=\vmatrix y & \tilde y \\ y''' & \tilde y''' \endvmatrix,
\\
u_4=\vmatrix y' & \tilde y' \\ y'' & \tilde y'' \endvmatrix,
\qquad
u_5=\vmatrix y' & \tilde y' \\ y''' & \tilde y''' \endvmatrix,
\qquad
u_6=\vmatrix y'' & \tilde y'' \\ y''' & \tilde y''' \endvmatrix.
\endgather
$$
Using the differential equation~\thetag{2.1}
for $y,\tilde y$, we get the system
$$
\gather
Y''=u_3+u_4, \qquad
u_3'=u_5-a_1Y-a_2Y'-a_3u_3,
\\
u_4'=u_5, \qquad
u_5'=u_6+a_0Y-a_2u_4-a_3u_5, \qquad
\\
u_6'=a_0Y'+a_1u_4-a_3u_6.
\tag3.1
\endgather
$$
Then
$$
\align
Y'''
&=u_3'+u_4'
=u_5-a_1Y-a_2Y'-a_3u_3+u_5
\\
&=2u_5-a_1Y-a_2Y'-a_3(Y''-u_4),
\tag3.2
\endalign
$$
and in the notation
$$
U=Y'''+a_3Y''+a_2Y'+a_1Y
$$
we may write \thetag{3.2} in the form
$$
2u_5=U-a_3u_4.
\tag3.3
$$
Therefore
$$
U'-a_3'u_4-a_3u_5=(U-a_3u_4)'=2u_5'
=2u_6+2a_0Y-2a_2u_4-2a_3u_5.
$$
Hence by \thetag{3.3}
$$
U'+\frac12a_3U-2a_0Y
=2u_6+\Bigl(\frac12a_3^2-2a_2+a_3'\Bigr)u_4
$$
and after taking the derivative
$$
U''+\frac12a_3U'+\frac12a_3'U-2a_0Y'-2a_0'Y
=2u_6'+\Bigl(\frac12a_3^2-2a_2+a_3'\Bigr)u_5
+(a_3a_3'-2a_2'+a_3'')u_4.
$$
Using \thetag{3.3} and \thetag{3.1} we write the last equality as
$W=Vu_4$, where
$$
W=U''+\frac32a_3U'+\Bigl(a_3+\frac12a_3^2-\frac14a_3^3\Bigr)U
-4a_0Y'-2(a_0'+a_0a_3)Y
$$
and
$$
V=a_3''-2a_2'+\frac32a_3a_3'-a_2a_3+\frac14a_3^2+2a_1.
$$
The condition $V=0$ is exactly \thetag{2.2}; if it holds,
we have $W=Vu_4=0$, which is the required 5th-order differential equation.
If $V\ne0$, then differentiate $W=Vu_4$ to get
$$
W'=V'u_4+Vu_5
=V'u_4+\frac12V(U-a_3u_4),
$$
hence
$$
VW'=\frac12V^2U+\Bigl(V'-\frac12a_3V\Bigr)W.
\tag3.4
$$
Equality \thetag{3.4} is the 6th-order linear differential equation for the
function $Y$; the coefficient of $Y^{(6)}$ in this equation
is $V\ne0$. \qed
\enddemo

From now on, we assume condition~\thetag{2.2} for a given 4th-order
MUM differential equation~\thetag{2.1}. To construct
a 5th-order MUM differential equation we modify
the above construction by taking
$$
w(z)=zY=z\vmatrix y & \tilde y \\ y' & \tilde y' \endvmatrix.
$$
Denote the resulting 5th-order differential equation for the function~$w$:
$$
w^{(5)}+b_4(z)w^{(4)}+b_3(z)w^{(3)}+b_2(z)w''+b_1(z)w'+b_0(z)w=0.
\tag3.5
$$

\remark{Remark}
Clearly, condition \thetag{2.2} induces a certain condition
for the coefficients of the differential equation~\thetag{3.5}.
We leave the corresponding (rather boring) relation to the reader
as a {\tt Maple} exercise.
\endremark

\proclaim{Proposition 4}
The differential equation \thetag{3.5} is MUM; its Frobenius
basis is given by the following formulas:
$$
\gathered
w_0=z\vmatrix y_0 & y_1 \\ y_0' & y_1' \endvmatrix,
\quad
w_1=z\vmatrix y_0 & y_2 \\ y_0' & y_2' \endvmatrix,
\quad
w_2=z\vmatrix y_0 & y_3 \\ y_0' & y_3' \endvmatrix
=z\vmatrix y_1 & y_2 \\ y_1' & y_2' \endvmatrix,
\\
w_3=\frac z2\vmatrix y_1 & y_3 \\ y_1' & y_3' \endvmatrix,
\quad
w_4=\frac z2\vmatrix y_2 & y_3 \\ y_2' & y_3' \endvmatrix.
\endgathered
\tag3.6
$$
\endproclaim

\remark{Remark}
The inverse is not true in the following sense. If we take
a MUM linear differential equation of order~5 with monodromy
group $O_5$, then its 4th-order differential pullback is not
necessarily MUM. An example to this is the 5th-order differential
operator
$$
D=\theta^5-6^6z\bigl(\theta+\tfrac16\bigr)\bigl(\theta+\tfrac26\bigr)
\bigl(\theta+\tfrac36\bigr)\bigl(\theta+\tfrac46\bigr)
\bigl(\theta+\tfrac56\bigr)
$$
(which is a natural generalization of the 4th-order operator
$$
D=\theta^4-5^5z\bigl(\theta+\tfrac15\bigr)\bigl(\theta+\tfrac25\bigr)
\bigl(\theta+\tfrac35\bigr)\bigl(\theta+\tfrac45\bigr)
$$
corresponding to the Calabi--Yau differential equation).
\endremark

\demo{Proof}
The five functions in~\thetag{3.6} are linearly independent solutions
to~\thetag{3.5}. Developing their $(\log z)$-expansions (with the
known structure of the Frobenius basis $y_0,y_1,y_2,y_3$ to~\thetag{2.1})
show that they form the Frobenius basis of \thetag{3.5}, hence the latter
differential equation is MUM.
\qed
\enddemo

\proclaim{Corollary \rm(Beukers's relations)}
We have
$$
2w_0w_4-2w_1w_3+w_2^2=0, \qquad
2w_0'w_4'-2w_1'w_3'+{w_2'}^2=0.
\tag3.7
$$
\endproclaim

\demo{Proof}
Write
$$
w_2^2=z\vmatrix y_0 & y_3 \\ y_0' & y_3' \endvmatrix
\cdot z\vmatrix y_1 & y_2 \\ y_1' & y_2' \endvmatrix.
$$
Then the first identity is trivial by expanding. In the same vein
we obtain the identity for the functions $\tilde w_j$, in which
the second line of the determinant contains $y_i'$s replaced by $y_i''$s.
Since $w_j'=w_j/z+\tilde w_j$, we finally arrive at the second identity
in~\thetag{3.7}.
\qed
\enddemo

\remark{Remark}
Beukers proved the above relations for a concrete
example of the 5th-order differential equation, which
we discuss in the next section.
\endremark

In the case of the 5th-order MUM differential equation~\thetag{3.5} we can
define analogues of a mirror map and a Yukawa coupling.
Namely, we can define $\tilde z(q)$ as the inverse of
$q(z)=\exp(w_1/w_0)$ and
$$
\tilde K(q)=\tilde N_0\cdot\frac{\d^2}{\d t^2}(w_2/w_0),
\qquad t=w_1/w_0.
$$
Developing the last $q$-series as the Lambert series,
$$
\tilde K(q)=\tilde N_0+\sum_{l=1}^\infty\frac{l^2\tilde N_lq^l}{1-q^l}
\tag3.8
$$
(we put $l^2$ instead of $l^3$ here),
we have the following experimental observation:
$\tilde z(q)\in\allowmathbreak K[[q]]$ and $\tilde N_l\in K$,
where $K$ is either $\Bbb Z$ or $\Bbb Z[1/p]$
for a certain prime~$p$ (the possibilities $p=3,5,7,23$ occur).

\head
4. Strange integrality related to $\zeta(4)$
\endhead

In the work \cite{Zu1}, the following 3-term polynomial
recursion is given:
$$
\align
&
(n+1)^5A_{n+1}
-3(2n+1)(3n^2+3n+1)(15n^2+15n+4)A_n
\\ &\qquad
-3n^3(3n-1)(3n+1)A_{n-1}=0
\qquad\text{for}\quad n\ge1.
\endalign
$$
If we take two linearly independent solutions $\{A_n\}$ and $\{B_n\}$
to the above recursion given by the initial data
$$
A_0=1, \quad A_1=12, \qquad\text{and}\qquad
B_0=0, \quad B_1=13,
$$
then
$$
\lim_{n\to\infty}\frac{B_n}{A_n}=\zeta(4)=\frac{\pi^4}{90}
$$
(see \cite{Zu1} for details). This recursion was previously
proved in \cite{Co} and \cite{So} but without any indication
of arithmetic properties of the sequences $\{A_n\}$ and $\{B_n\}$.
In~\cite{Zu1}, it is proved that
$$
D_nA_n,D_n^5B_n\in\Bbb Z,
\qquad\text{where}\quad
D_n=\text{the least common multiple of}\ 1,2,\dots,n,
$$
and conjectured some stronger inclusions that were finally proved
in \cite{KR}. In particular, from \cite{KR} we have the inclusions
$$
A_n,D_n^4B_n\in\Bbb Z,
$$
and even the following explicit formula (see also~\cite{Zu4}):
$$
A_n
=\sum_{j,k}\binom nj^2\binom nk^2\binom{n+j}n\binom{n+k}n\binom{j+k}n,
$$
where the binomial coefficients $\binom ab$ are zero if
$b<0$ or $a<b$.

A natural object related to the above recursion
is the 5th-order MUM differential equation
$$
\align
&
z^5(27z^2+270z-1)w^{(5)}
+z^4(405z^2+3375z-10)w^{(4)}
\\ &\qquad
+z^3(1752z^2+11502z-25)w^{(3)}
+z^2(2412z^2+11259z-15)w''
\\ &\qquad
+z(816z^2+2130z-1)w'
+12z(2z+1)w=0.
\tag4.1
\endalign
$$
Denote by $w_0,w_1,w_2,w_3,w_4$ the Frobenius basis of \thetag{4.1}
near $z=0$, so that
$$
w_0
=w_0(z)=\sum_{n=0}^\infty A_nz^n
=1+12z+804z^2+88680z^3+12386340z^4+\dotsb.
$$
It was pointed out by F.~Beukers~\cite{Be4} that
the (Zariski closure of the) monodromy group
of the differential equation~\thetag{4.1}
turns out to be $O_5$. This result (as well as identities \thetag{3.7}
for the latter functions $w_j$s) are consequences of the following
statement.

\proclaim{Proposition 5}
The $5$th-order MUM differential equation~\thetag{4.1}
is constructed from the $4$th-order MUM differential
equation~\thetag{2.1} with
$$
\align
a_3
&=\frac{6(36z^2+315z-1)}{z(27z^2+270z-1)},
\\
a_2
&=\frac{10530z^4+180387z^3+759417z^2-4671z+7}{z^2(27z^2+270z-1)^2},
\\
a_1
&=\frac{96228z^6+2421009z^5+18416565z^4+49339854z^3
-450054z^2+1107z-1}
{z^3(27z^2+270z-1)^3},
\\
a_0
&=\frac{\aligned
3(8748z^7+253692z^6+10875303z^5-37601010z^4
\qquad \\ \vspace{-4.5pt}
+13643328z^3+432135z^2+13224z+11)
\endaligned}
{z^3(27z^2+270z-1)^4}
\endalign
$$
by the algorithm of Section~\rom3.
\endproclaim

\demo{Proof}
By direct computation (using {\tt Maple}).
\qed
\enddemo

\remark{Remark}
The original proof by Beukers~\cite{Be5} of relations \thetag{3.7}
for the differential equation \thetag{4.1} follows the
principle `the proof is trivial once has found the explicit formula'.
We reproduce the arguments from \cite{Be5} here.
The symmetric square of the 5th order equation \thetag{4.1}
has order 14, one order less than expected~15. This means
that there should be a homogeneous quadratic relation
between the solutions $w_0,\dots,w_4$ of~\thetag{4.1}.
It suffices to find a relation which vanishes at $z=0$
of order one higher than the local exponents given by the
symmetric square of order~14. (Clearly, any non-trivial
relation satisfies the 14th order equation; if we find one with
one vanishing order higher than predicted by this 14th order equation,
it must be the trivial solution.) If we replace the $w_j$s by
their derivatives, we get quadratic forms that form a one-dimensional
monodromy representation. Looking at the local exponents
we can infer that they should be rational functions
with poles only at infinity. Hence the expressions are
polynomial of degree again given by the local exponents
at infinity of the symmetric square.
\endremark

\medskip
Here are some observations related to the $5$th-order
differential equation~\thetag{4.1}
and its $4$th-order pullback given in Proposition~5 (we denote
by $y_0,y_1,y_2,y_3$ the Frobenius basis of the latter
differential equation).

\example{Observation 1}
We have $y_0\in\Bbb Z[1/2][[z]]$ but not $y_0\in\Bbb Z[[z]]$,
while $z(q),K(q)\in\Bbb Z[[q]]$ and for $K(q)$ we have
the Lambert expansion~\thetag{2.9} with $N_l\in\Bbb Z$
but $N_l/N_0<0$.
\endexample

\example{Observation 2}
We have $w_0=u_0^2$, where
$$
u_0=1+6z+384z^2+42036z^3+5867226z^4+\dotsb
\in\Bbb Z[[z]].
$$
\endexample

\example{Observation 3}
One has the following generalized Beukers's relations:
for $k\ge2$
$$
2w_0^{(k)}w_4^{(k)}-2w_1^{(k)}w_3^{(k)}+(w_2^{(k)})^2
=\sum_\nu c_{k,\nu}(9z)^\nu,
$$
where $c_{k,\nu}\in\Bbb Z$.
\endexample

\example{Observation 4}
Put $t=w_1/w_0$, $q(z)=e^t=z+O(z^2)$ and denote by $\tilde z(q)$ its inverse.
Then $q(z)\in\Bbb Z[[z]]$
and, hence, $\tilde z(q)\in\Bbb Z[[q]]$. Moreover,
$$
\align
q(z)&=z(1+15z+1145z^2+\dotsb)^3,
\\
\tilde z(q)&=q(1-15q-245q^2-20138q^3-2043703q^4+\dotsb)^3.
\endalign
$$
\endexample

\example{Observation 5}
Expand $w_0$ in Lambert series:
$$
2w_0(z)
=2+\sum_{n=1}^\infty\frac{n^3L_nz^n}{1-z^n}.
$$
Then $L_n\in\Bbb Z$ and $L_n>0$.
What is a enumerative meaning of the
numbers $L_n$ (i.e., do they count some geometric objects)?
\endexample

\example{Observation 6}
Write
$$
\tilde K(q)=\frac{\d^2}{\d t^2}(w_2/w_0)
=1+\sum_{l=1}^\infty\frac{l^2\tilde N_lq^l}{1-q^l}.
$$
Then $\tilde N_l\in\Bbb Z$ and $\tilde N_l>0$.
What is a enumerative meaning of the
numbers $\tilde N_l$ (i.e., do they count some geometric objects)?
\endexample

\example{Observation 7 \rm(Kummer supercongruences)}
If we write $w_0$ and $\tilde K$ as functions of~$z$,
$$
w_0(z)=\sum_{n=0}^\infty A(n)z^n,
\qquad
\tilde K(q(z))=\sum_{n=0}^\infty C(n)z^n,
$$
then
$$
A(np^r)\equiv A(np^{r-1})\pmod{p^{3r}}, \qquad
C(np^r)\equiv C(np^{r-1})\pmod{p^{2r}},
$$
for primes $p$ and $(n,p)=1$; $r=1,2,\dots$\,.
\endexample

Recently the second author~\cite{Zu5} (see also Theorem~4
in~\cite{Zu2}), when approximating
$\Li_3(z)$ (see definition~\thetag{9.3} below), found
two more {\it fourth order\/} examples related
to simultaneous rational approximations
to $\zeta(2)$ and~$\zeta(3)$. About the same time the first
author discovered compact formulas
$$
A_n=\sum_{k,l}\binom nk^2\binom nl^2\binom{k+l}l\binom{n+k}n
\qquad n=0,1,2,\dots,
\tag4.2
$$
and
$$
A_n=\binom{2n}n\sum_k\binom nk^2\binom{n+k}n\binom{n+2k}n,
\qquad n=0,1,2,\dots,
$$
for generating series $y_0(z)=\sum_{n=0}^\infty A_nz^n$
satisfying exactly the same differential equations.
These examples are cases \#195 and \#209 in \cite{AESZ, Table~A}.
The coefficients~\thetag{4.2} is of the form similar
to those in~\cite{BS2} (cf\. case \#26),
and van Straten has confirmed
that equation \#195 indeed `comes from geometry'.
Thus we have a direct link between number theory and geometry.

\head
5. Strategy of finding new Calabi--Yau equations
\endhead

How much do we know about Calabi--Yau differential equations?
We discovered a list of 29 interesting cases taken from different
papers on the subject; these cases initiated
our systematization of Calabi--Yau equations and started
the table in~\cite{AESZ}%
\footnote{Unfortunately, we are not
quite sure that the nice arithmetic properties of mirror maps
and Yukawa couplings are really proved in all presented cases!
In this sense, \cite{AESZ, Table~A} involves cases that are only expected
to be `arithmetically nice'.}. The first 14 cases
in the table correspond to all known hypergeometric
Calabi--Yau equations, while the origin of the remaining 15 known
cases is due to certain variations of multiple hypergeometric
series (or multiple binomial sums). The latter circumstance
prompted us to look on those linear MUM differential equations,
whose analytic solutions $y_0(z)\in\Bbb Z[[z]]$
are single or multiple sums of certain hypergeometric terms.
There is an algorithmic way, which is due to Gosper and Zeilberger
known as {\it creative telescoping}, to find
polynomial recursions for such sums \cite{PWZ}. Unfortunately,
the algorithm for multiple summations (being a part
of the general {\it WZ-theory\/} developed
by H.~Wilf and D.~Zeilberger~\cite{WZ}) is not implemented
in {\tt Maple} but solving auxiliary systems of linear equations is.
The recursions may be easily transfered to linear differential
equations satisfied by $y_0(z)$.
That was our general recipe for verifying known cases \#15--29,
and then we continued our experiments by constructing several
new cases, also in combination with other hypergeometric
methods---quadratic transformations (Section~6),
Hadamard products (Section~7) and other transformations (Section~8).
In several cases the algorithm
of creative telescoping led to reducible differential equations,
but then we used {\tt DEtools\/} implemented in {\tt Maple}
to factorize corresponding differential operators.

The way {\it a l\`a Ap\'ery\/} of producing new examples
of the mirror maps and Yukawa couplings, decribed in the previous
section seems to be very occasional: we do not really have
other examples of rational approximations to `interesting' constants
(but see Section~10). Nevertheless, the method of finding
Calabi--Yau equations by pulling back some `nice' 5th-order
linear differential equations works in several other situations
as well. We stress that the 5th- to 4th-order reduction
did not appear before in connection with problems of
Calabi--Yau manifolds.

Our table in~\cite{AESZ} extends the list of 29 known cases
and contains a sporadic list of more than two hundred Calabi--Yau
(conjecturally) differential equations.
The table is long but we have found many more equations
that we decided not to include (not to stretch the reader's
patience too far!). In particular, we have many examples of type
$A_n=\sum_k(-1)^k2^{n-2k}\binom n{2k}\binom{2k}k\bigl(\,\cdot\,\bigr)$.
The most unfortunate thing seems to be the fact that
we cannot see any reasonable end of our list
and, therefore, we do not expect to be much criticized
by the reader for missing a hundred or more other interesting cases.

The analysis of \cite{AESZ, Table~A} shows that, besides the 14
hypergeometric functions (cases \#1--14), we have a plenty
of other `simple' Calabi--Yau differential equations of type
$$
\theta^4-c_1z(a\theta^2+a\theta+b)P(\theta)+c_2z^2Q(\theta),
$$
where $P(\theta)$ is of the form $(2\theta+1)^2$,
$(3\theta+1)(3\theta+2)$, $(4\theta+1)(4\theta+3)$ or
$(6\theta+1)(6\theta+5)$, and $Q(\theta)$ splits
into linear factors of similar shapes.
Here $a\in\{2,3,5,7,8,9,11,13,17\}$ and $c_1,c_2$ are usually
powers of~$2$ and (or)~$3$.
We also have $5$th-order differential equations
similar to the the $4$th-order in~2. Their 4th-order
pullbacks (derived by the method of Section~3)
are much more complicated. Case \#32 seems to be even
more complicated:
$$
\align
&
\theta^5
-3z(2\theta+1)(3\theta^2+3\theta+1)(15\theta^2+15\theta+4)
\\ &\qquad
-3z^2(\theta+1)^3(3\theta+2)(3\theta+4)
\endalign
$$
(this is the differential equation that started this paper).
Then there are the cases, when both the 4th-order differential
equation and the 5th-order equation (obtained by the wronskian
construction) are very complicated. The corresponding
linear differential operators usually appear as factors
of higher order differential operators, for which we have
no pattern at all.

One can suspect that the equations in the list from
\cite{AESZ, Table~A} are Picard--Fuchs equations of one-dimensional
families of Calabi--Yau threefolds. In many cases, such as
the hypergeometric ones and the Hadamard products, this should
not be hard to see. For other cases, such as the one related
to~$\zeta(4)$, this fact will be quite non-trivial, but perhaps
all the more interesting. This is a reason of why the equations
are named Calabi--Yau equations.

\remark{Remark}
We are quite surprised by the enomous amount of
the 4th-order linear differential equations with
(expectively) nice arithmetic properties.
After discovering equation \#9, the first author
for a long time thought that there were no more
examples but \#1--14. Then after finding the papers
\cite{BS1}, \cite{BS2} and examples \#15--28 there,
it was still hard to find any new cases. The use
of the Gosper--Zeilberger algorithm changed this
drastically. For a while we thought that whenever
$A_n$ was a binomial sum (even multiple) and it gave
a 4th-order MUM differential equation, then we
automatically got an integer mirror map $z(q)$
and the Yukawa coupling $K(q)$ as in~\thetag{2.9}
with integers $N_l$. There is, however, the following
example: the generating series $y_0=\sum_{n=0}^\infty A_nz^n$, where
$$
A_n=\sum_k(-1)^k4^{n-4k}\binom n{4k}\frac{(5k)!}{k!^5},
$$
satisfies a 4th-order MUM differential equation with polynomial
coefficients of degree~8; the equation produces the non-integer
mirror map~$z(q)$ (actually, $z(q)\in\allowmathbreak\Bbb Z[1/2][[q]]$)
but the Yukawa coupling has integers $N_l$. There are more
examples, also for the cases of non-integers $N_l$; we indicate
one in Section~7 below.
\endremark

\head
6. Quadratic transformations and mirror maps
\endhead

A purely analytic, hypergeometric
machinery---quadratic and higher-order transfor\-mations---may
be put forward for constructing new examples from the old
ones. The method powerfully works in the case of 2nd- and 3rd-order
differential equations, when we have to substitute a suitable
modular function in place of variable~$z$; see, for example, \cite{HM}.

Our starting point is the existence
of quadratic transformations for ${}_2F_1$- and ${}_3F_2$-series,
like
$$
{}_2F_1\biggl(\matrix\format&\,\c\\
a, & b \\ & 1+a-b
\endmatrix\biggm|z\biggr)
=(1-z)^{-a}\cdot{}_2F_1\biggl(\matrix\format&\,\c\\
\frac12a, & \frac12+\frac12a-b \\ & 1+a-b
\endmatrix\biggm|-\frac{4z}{(1-z)^2}\biggr)
\tag6.1
$$
due to Gau\ss\ and
$$
\alignat1
&
{}_3F_2\biggl(\matrix\format&\,\c\\
a, & b, & c \\ & 1+a-b, & 1+a-c
\endmatrix\biggm|z\biggr)
\\ &\qquad
=(1-z)^{-a}\cdot{}_3F_2\biggl(\matrix\format&\,\c\\
\frac12a, & \frac12+\frac12a, & 1+a-b-c \\ & 1+a-b, & 1+a-c
\endmatrix\biggm|-\frac{4z}{(1-z)^2}\biggr)
\tag6.2
\endalignat
$$
due to Whipple. Substituting $a=b=c=\frac12$ in them, we get algebraic
relations between different mirrors that are necessarily {\it modular\/}
maps in these cases. These quadratic transformations have a natural
generalization to the case of higher dimensional hypergeometric series
(see \cite{Zu3}, Sections~6 and~7). Similar but particular results
are given in the following statement.

\proclaim{Proposition 6}
The following quadratic transformations are available
for ${}_4F_3$-series:
$$
\align
&
{}_4F_3\biggl(\matrix
a, \, b, \, c, \, d \\ 1+a-b, \, 1+a-c, \, 1+a-d \endmatrix
\biggm|z\biggr)
\\ &\qquad
=\frac1{(1+z)^a}
\sum_{n=0}^\infty\frac{(\frac12a)_n\,(\frac12+\frac12a)_n}
{(1+a-b)_n\,(1+a-c)_n}
\biggl(\frac{4z}{(1+z)^2}\biggr)^n
\\ &\qquad\quad\times
\sum_{\nu=0}^n\frac{(b)_\nu\,(c)_\nu\,(1+a-b-c)_{n-\nu}}
{\nu!\,(n-\nu)!\,(1+a-d)_\nu}
\tag6.3
\\ \allowdisplaybreak &\qquad
=\frac1{(1-z)^a}
\sum_{n=0}^\infty\frac{(\frac12a)_n}{(1+a-b)_n}
\biggl(-\frac{4z}{(1-z)^2}\biggr)^n
\sum_{\mu=0}^n
\frac{(b)_\mu\,(\frac12+\frac12a)_\mu\,(\frac12+\frac12a-b)_{n-\mu}}
{(n-\mu)!\,(1+a-c)_\mu}
\\ &\qquad\quad\times
\sum_{\nu=0}^\mu\frac{(c)_\nu}{\nu!\,(\mu-\nu)!\,(1+a-d)_\nu}(-1)^\nu.
\tag6.4
\endalign
$$
\endproclaim

\demo{Proof}
Writing
$$
\align
&
{}_4F_3\biggl(\matrix
a, \, b, \, c, \, d \\ 1+a-b, \, 1+a-c, \, 1+a-d \endmatrix
\biggm|z\biggr)
\\ &\qquad
=\sum_{n=0}^\infty\frac{(a)_n\,(b)_n\,(c)_n}
{n!\,(1+a-b)_n\,(1+a-c)_n}(-z)^n\cdot{}_2F_1\biggl(\matrix
-n, \ a+n \\ 1+a-d \endmatrix\bigg|1\biggr)
\\ &\qquad
=\sum_{\nu=0}^\infty\frac{(-1)^\nu}{\nu!\,(1+a-d)_\nu}
\sum_{n=\nu}^\infty\frac{(a)_{\nu+n}\,(b)_n\,(c)_n}
{(n-\nu)!\,(1+a-b)_n\,(1+a-c)_n}(-z)^n
\\ &\qquad
=\sum_{\nu=0}^\infty\frac{(a)_{2\nu}(b)_\nu(c)_\nu}
{\nu!\,(1+a-b)_\nu\,(1+a-c)_\nu\,(1+a-d)_\nu}
z^\nu
\\ &\qquad\quad\times
{}_3F_2\biggl(\matrix
a+2\nu, \, b+\nu, \, c+\nu \\ 1+a-b+\nu, \, 1+a-c+\nu \endmatrix
\biggm|-z\biggr).
\tag6.5
\endalign
$$
To the latter ${}_3F_2$-series we apply the quadratic
transformation~\thetag{6.2} and reorder summations
to get~\thetag{6.3}.

For the proof of~\thetag{6.4}, we proceed as in~\thetag{6.5} to deduce
$$
\alignat1
&
{}_3F_2\biggl(\matrix
a, \, b, \, c \\ 1+a-b, \, 1+a-c \endmatrix
\biggm|z\biggr)
\\ &\qquad
=\sum_{\lambda=0}^\infty\frac{(a)_{2\lambda}(b)_\lambda}
{\lambda!\,(1+a-b)_\lambda\,(1+a-c)_\lambda}
z^\lambda
\cdot{}_2F_1\biggl(\matrix
a+2\lambda, \, b+\lambda \\ 1+a-b+\lambda \endmatrix
\biggm|-z\biggr).
\tag6.6
\endalignat
$$
To the inner ${}_2F_1$-series we now apply the
transformation~\thetag{6.1} to get, as before, the double
series for the left-hand side of~\thetag{6.6}. Then substituting
the resulting formula for the ${}_3F_2$-series into~\thetag{6.5}
we arrive at the desired identity~\thetag{6.4}.
\qed
\enddemo

Plugging in $a=b=c=d=\frac12$ we obtain
$$
\align
\sum_{n=0}^\infty\binom{2n}n^4\biggl(\frac z{2^8}\biggr)^n
&=(1+z)^{-1/2}
\sum_{n=0}^\infty A_n^{(+)}\biggl(\frac z{2^8(1+z)^2}\biggr)^n
\\
&=(1-z)^{-1/2}
\sum_{n=0}^\infty A_n^{(-)}\biggl(-\frac z{2^8(1-z)^2}\biggr)^n,
\endalign
$$
where
$$
\align
A_n^{(+)}
&=\frac{(4n)!}{(2n)!\,n!^2}
\sum_{\nu=0}^n2^{2(n-\nu)}\binom{2\nu}\nu^2\binom{2(n-\nu)}{n-\nu},
\tag6.7
\\
A_n^{(-)}
&=2^{2n}\frac{\prod_{j=0}^{n-1}(1+4j)}{n!}
\sum_{\mu=0}^n2^{4(n-\mu)}\binom{2\mu}\mu
\frac{\prod_{j=0}^{n-\mu-1}(1+4j)}{(n-\mu)!}
\,\frac{\prod_{j=0}^{\mu-1}(3+4j)}{\mu!}
\\ &\qquad\times
\sum_{\nu=0}^\mu2^{2(\mu-\nu)}\binom{2\nu}\nu\binom\mu\nu(-1)^\nu.
\tag6.8
\endalign
$$

The resulting series $\sum_{n=0}^\infty A_n^{(+)}z^n$
and $\sum_{n=0}^\infty A_n^{(-)}z^n$ are solutions of
4th-order linear differential equations. They gave
two new examples (cases \#30 and \#31)
of differential equations producing
mirror maps and Yukawa couplings with desired integrality
properties. It seems to be an interesting problem to find
quadratic and high-order transformations for other cases
from \cite{AESZ, Table~A}.

\medskip
Besides hypergeometric transformations, there is also
an `obvious' way to produce infinitely many linear
MUM differential equations starting with the only one such
example~\thetag{1.1}.

\proclaim{Proposition 7}
Let $y=y(z)$ be a generic solution of equation~\thetag{1.1}
of order~$s$ with coefficients in $\Bbb C(z)$
and let $u=u(z)$ be any function satisfying $u'/u\in\Bbb C(z)$.
Then the function $w=y/u$ is a solution of the linear
differential equation of order~$s$ with rational coefficients,
and this equation depends only on \thetag{1.1} but not
on its solution~$y$.
\endproclaim

\demo{Proof}
Substituting $y=uw$ in \thetag{1.1} we get the equation
$$
(uw)^{(s)}+a_{s-1}(z)(uw)^{(s-1)}+\dots+a_1(z)(uw)'+a_0(z)uw=0
\tag6.9
$$
in $w$ unknown. Developing the derivatives by the formulas
$$
(uw)^{(k)}=\sum_{j=0}^k\binom kju^{(j)}w^{(k-j)},
\qquad k=1,2,\dots,s,
$$
and dividing the left-hand side in~\thetag{6.9} by $u$,
we obtain the $s$th-order linear differential equation for~$w$
with rational coefficients since $u^{(j)}/u\in\Bbb C(z)$
for any $j=0,1,2,\dots$\,.
\qed
\enddemo

By Proposition~7, taking a known example of a Calabi--Yau equation
with the corresponding series $y_0=\sum_{n=0}^\infty A_nz^n\in\Bbb Z[[z]]$
and choosing any $u\in1+z\Bbb Z[[z]]$ with the property $u'/u\in\Bbb Q(z)$,
we see that $w_0=y_0/u$ is in $\Bbb Z[[z]]$
and $w=y/u$ also satisfies a Calabi--Yau equation, with the same
mirror map $z(q)$ and Yukawa coupling $K(q)$, since
$w_j/w_0=y_j/y_0$ for $j=1,2,3$.
For instance, the choice $u=\sqrt{1-4pz}$, where $p$ is an arbitrary
integer, leads to an infinite family of Calabi--Yau equations
for the sequences
$$
\wh A_n=\sum_{k=0}^np^{n-k}\binom{2n-2k}{n-k}A_k.
$$
Clearly, examples of such type are not presented in Appendix~A.

\medskip
We also found experimentally another (doubly) infinite
classes of Calabi--Yau differential equations depending
on positive integer parameters $p$ and~$r$.

\proclaim{Proposition 8}
Given $A_n$, one of the cases presented in \cite{AESZ, Table~A} with
the corresponding Yukawa coupling~$K(q)$, define
$$
\wh A_n=\sum_kp^{n-rk}\binom n{rk}A_k.
$$
Then $\wh y_0=\sum_{n=0}^\infty\wh A_nz^n$ satisfies a $4$th-order
MUM differential equation with Yukawa coupling $\wh K(q)=K(q^r)$.
\endproclaim

\remark{Remark}
It can occur that the mirror map $\wh z(q)$ does not lie
in $\Bbb Z[[q]]$ (but the denominators of its coefficients
are necessarily divisors of powers of~$r$).
\endremark

\demo{Proof}
We only give a sketch of the proof, without indicating computational details.
Suppose that $D$~is the linear 4th-order differential operator
annihilating the given series $y_0(z)=\sum_{n=0}^\infty A_nz^n$
and $y_0,y_1,y_2,y_3$ is the Frobenius basis of the differential
equation $Dy=0$. Change the variable
$z\mapsto Z(z)=\bigl(z/(1-pz)\bigr)^r$; then substituting
$$
\delta=z\frac{\d}{\d z}=\frac r{1-pz}\cdot Z\frac{\d}{\d Z}
$$
into the differential equation $Dy=0$ leads to a new
4th-order MUM differential equation that annihilates
$Y_0(z)=y_0(Z(z))$. There is no difficulty in computing
the Frobenius basis of the new differential equation
(cf.~\thetag{1.2}): it is
$$
Y_j(z)=r^{-j}y_j(Z(z)), \qquad j=0,1,2,3.
$$
Therefore, $T(z)=Y_1(z)/Y_0(z)=\frac1rt(Z(z))$ and the new mirror
map $Z(q)$ is related to the old one $z(q)$ in accordance with the formula
$$
z(q^r)=\biggl(\frac{Z(q)}{1-pZ(q)}\biggr)^r.
$$
Finally,
$$
\frac{\d^2}{\d T^2}\biggl(\frac{Y_2}{Y_0}\biggr)=K(q^r).
$$
Applying now the composition
$$
y(z)\mapsto\wt y(z)=\frac1{1-pz}\cdot y(Z(z))
=\frac1{1-pz}\cdot y\biggl(\biggl(\frac z{1-pz}\biggr)^r\biggr)
$$
with the help of Proposition~7 we conclude that
the differential equation for~$\wt y(z)$ is of order~4 and MUM,
and the corresponding Yukawa coupling is $K(q^r)$. It remains
to verify that
$$
\wt y_0(z)
=\sum_{n=0}^\infty\wt A_nz^n
=\frac1{1-pz}\sum_{n=0}^\infty A_n\biggl(\frac z{1-pz}\biggr)^{rn}
=\frac1{1-pz}\cdot y_0\biggl(\biggl(\frac z{1-pz}\biggr)^r\biggr).
\qed
$$
\enddemo

\remark{Remark}
As pointed out to us by Beukers, Propositions~7 and 8
might be discarded by choosing a normalization
for Calabi--Yau differential equations. Natural requirements are
as follows:
\roster
\item"(a)" the smallest local exponent at the finite
singularities is zero;
\item"(b)" the point $z=\infty$ is a singularity.
\endroster
However, the propositions remain useful if these conditions
are abandoned.
\endremark

\head
7. Hadamard products
\endhead

Using the algorithm of creative telescoping we found
almost a hundred new Calabi--Yau differential equations.
Then D.~van Straten suggested that using Hadamard products
of solutions to `nice' 2nd-order equations could give required
4th-order equations. Surprisingly, this was also the case
for many of the examples we had found just by accident.

Let
$$
u=\sum_{n=0}^\infty b_nz^n, \qquad v=\sum_{n=0}^\infty c_nz^n
$$
be two $D$-finite (i.e., satisfying a linear differential
equation of finite order with polynomial coefficients)
power series. Then the Hadamard product
$$
y=\sum_{n=0}^\infty a_nz^n=u*v=\sum_{n=0}^\infty b_nc_nz^n
\tag7.1
$$
is also $D$-finite (see \cite{Sta, p.~194}). If $D_u$ and $D_v$
are linear differential operators annihilating $u$ and $v$, respectively,
by $D_u*D_v$ we denote the differential operator annihilating \thetag{7.1}.
We do not know a general algorithm for computing $D_u*D_v$,
but for a given pair of operators $D_u,D_v$ the problem is
easily solved by linear algebra arguments.
We have found about 30 second order MUM differential equations
coming from binomial coefficients; Zagier's manuscript~\cite{Za}
provides us with 36 second order examples with `nice' arithmetic properties,
although a closed binomial-sum formula is not known in all cases.
Consider the following examples
for which we were successful in finding
the differential equations for their Hadamard products:
$$
\allowdisplaybreaks
\xxalignat2
&\text{(a)} &
A_n&=\sum_{k=0}^n\binom nk^3, \quad
&
D&=\theta^2-z(7\theta^2+7\theta+2)-8z^2(\theta+1)^2;
\\
&\text{(b)} &
A_n&=\sum_{k=0}^n\binom nk^2\binom{n+k}k, \quad
&
D&=\theta^2-z(11\theta^2+11\theta+3)-z^2(\theta+1)^2;
\\
&\text{(c)} &
A_n&=\sum_{k=0}^n\binom nk^2\binom{2k}k, \quad
&
D&=\theta^2-z(10\theta^2+10\theta+3)+9z^2(\theta+1)^2;
\\
&\text{(d)} &
A_n&=\sum_{k=0}^n\binom nk\binom{2k}k\binom{2n-2k}{n-k}, \quad
&
D&=\theta^2-4z(3\theta^2+3\theta+1)+32z^2(\theta+1)^2;
\\
&\text{(e)} &
A_n&=\sum_{k=0}^n4^{n-k}\binom{2k}k^2\binom{2n-2k}{n-k}, \quad
&
D&=\theta^2-4z(8\theta^2+8\theta+3)+256(\theta+1)^2.
\endxxalignat
$$
Taking all possible Hadamard products of these examples (including squares)
we get 15 new differential equations. Squares will be of degree~5 and the
others of degree~8. The products are listed in \#100--107, 113--115 and
120--123.

We may also extend the above list by the following examples
(f)--(h) due to D.~Zagier~\cite{Za} and (i), (j) due to
C.~van Enckevort and D.~van Straten~\cite{ES}:
$$
\allowdisplaybreaks
\xxalignat2
&\text{(f)} &
A_n&=\sum_{k=0}^{\[n/3\]}(-1)^k3^{n-3k}\binom n{3k}\frac{(3k)!}{k!^3}, \quad
D=\theta^2-3z(3\theta^2+3\theta+1)+27z^2(\theta+1)^2;
\\
&\text{(g)} &
& \kern52mm
D=\theta^2-z(17\theta^2+17\theta+6)+72z^2(\theta+1)^2;
\\
&\text{(h)} &
A_n&=27^n\sum_k(-1)^k\binom{-2/3}k\binom{-1/3}{n-k}^2,
\\ &&
D&=\theta^2-3z(18\theta^2+18\theta+7)+729z^2(\theta+1)^2;
\\
&\text{(i)} &
A_n&=64^n\sum_k(-1)^k\binom{-3/4}k\binom{-1/4}{n-k}^2,
\\ &&
D&=\theta^2-4z(32\theta^2+32\theta+13)+4096z^2(\theta+1)^2;
\\
&\text{(j)} &
A_n&=432^n\sum_k(-1)^k\binom{-5/6}k\binom{-1/6}{n-k}^2,
\\ &&
D&=\theta^2-12z(72\theta^2+72\theta+31)+186624z^2(\theta+1)^2.
\endxxalignat
$$
In case (g), no explicit formulas for $A_n$ is known.
In other cases the explicit formulas for $A_n$ were
found in the following way.

The Legendre function
$$
P_a(z)
={}_2F_1\biggl(\matrix
-a, & a+1 \\ & 1 \endmatrix
\biggm|\frac{1-z}2\biggr)
$$
satisfies the differential equation
$$
(1-z^2)\frac{\d^2y}{\d z^2}-2z\frac{\d y}{\d z}+a(a+1)y=0.
$$
Therefore, the function
$$
y_0(z)=\frac1{1-z}P_a\biggl(\frac{1+z}{1-z}\biggr)
$$
is annihilated by the differential operator
$$
\theta^2-z(2\theta^2+2\theta+a^2+a+1)+z^2(\theta+1)^2
$$
(see \cite{Za}). We have
$$
(1-z)^aP_a\biggl(\frac{1+z}{1-z}\biggr)
=\sum_{n=0}^\infty\binom an^2z^n
$$
which gives
$$
\frac1{1-cz}P_a\biggl(\frac{1+cz}{1-cz}\biggr)
=\sum_{n=0}^\infty A_nz^n,
$$
where
$$
A_n=c^n\sum_k(-1)^k\binom{-1-a}k\binom a{n-k}^2.
$$
Putting $a=-1/3,-1/4,-1/6$ for a suitable choice of $c$
we get cases (h)--(j).

Let further
$$
\xxalignat2
&\text{(A)} &
A_n&=\binom{2n}n^2, \quad
&
D&=\theta^2-4z(2\theta+1)^2;
\\
&\text{(B)} &
A_n&=\frac{(3n)!}{n!^3}, \quad
&
D&=\theta^2-3z(3\theta+1)(3\theta+2);
\\
&\text{(C)} &
A_n&=\frac{(4n)!}{n!^2(2n)!}, \quad
&
D&=\theta^2-4z(4\theta+1)(4\theta+3);
\\
&\text{(D)} &
A_n&=\frac{(6n)!}{n!(2n)!(3n)!}, \quad
&
D&=\theta^2-12z(6\theta+1)(6\theta+5).
\kern40mm
\endxxalignat
$$
One can prove the following formula
for the Hadamard product by the method presented below
($P$ and $Q$ are quadratic polynomials
and $c$ is an arbitrary constant):
$$
\bigl(\theta^2-zP(\theta)-cz^2(\theta+1)^2\bigr)
*\bigl(\theta^2-zQ(\theta)\bigr)
=\theta^4-zP(\theta)Q(\theta)-cz^2Q(\theta+1)Q(\theta).
$$
Then we get the following table of cases\footnote{%
We will not write out in~\cite{AESZ, Table~A}
the equations in cases \#133--143,
since the reader can easily derive them himself by the formula above.
We have, however, checked that one gets integers everywhere.}
from \cite{AESZ, Table~A} corresponding to Hadamard products:
\medskip
\line{\hss\vbox{\offinterlineskip
\halign to55mm{\strut\tabskip=100pt minus 100pt
\strut\vrule\vphantom{\vrule height10.1pt}#&\hbox to10.5mm{\hfil#\hfil}&%
\vrule#&\hbox to10.5mm{\hfil#\hfil}&%
\vrule#&\hbox to10.5mm{\hfil#\hfil}&%
\vrule#&\hbox to10.5mm{\hfil#\hfil}&%
\vrule#&\hbox to10.5mm{\hfil#\hfil}&%
\vrule#\tabskip=0pt\cr\noalign{\hrule}
&     && (A) && (B) &&         (C) && (D) &\cr
\noalign{\hrule\vskip1pt\hrule}
& (a) &&  45 &&  15 &&          68 &&  62 &\cr
& (b) &&  25 &&  24 &&          51 &&  63 &\cr
& (c) &&  58 &&  70 &&          69 &&  64 &\cr
& (d) &&  36 &&  48 &&          38 &&  65 &\cr
& (e) && 111 && 110 &&          30 && 112 &\cr
& (f) && 133 && 134 &&         135 && 136 &\cr
& (g) && 137 && 138 &&         139 && 140 &\cr
& (h) && 141 && 142 && $\emptyset$ && 143 &\cr
\noalign{\hrule}
}}\hss}
\medskip
\noindent
The empty set at place (C)$*$(h) means that the Yukawa
coupling $K(q)$ is constant.

In our computations we applied the following method
communicated to us by A.~Meurman. The method is, in a sense,
a discrete analogue of the wronskian formalism.

\remark{Method to find Hadamard products}
Assume
$$
A_{n+2}=P_1(n)A_{n+1}+Q_1(n)A_n, \quad
B_{n+2}=P_2(n)B_{n+1}+Q_2(n)B_n,
$$
and we want to find a recursion formula for $C_n=A_nB_n$.
For this, define
$$
\gather
R_j(n)=P_j(n+1)P_j(n)+Q_j(n+1), \quad
S_j(n)=P_j(n+1)Q_j(n), \qquad j=1,2,
\\
\gathered
U_j(n)=P_j(n+2)R_j(n)+Q_j(n+2)P_j(n), \\
V_j(n)=P_j(n+2)S_j(n)+Q_j(n+2)Q_j(n),
\endgathered
\qquad j=1,2,
\endgather
$$
and take
$$
\gather
T_2(n)=R_2S_1U_1V_2-R_1S_2U_2V_1, \quad
T_3(n)=P_1Q_2U_2V_1-P_2Q_1U_1V_2, \\
T_4(n)=R_1S_2P_2Q_1-R_2S_1P_1Q_2,
\endgather
$$
and
$$
W_0(n)=Q_1Q_2T_2+S_1S_2T_3+V_1V_2T_4, \quad
W_1(n)=P_1P_2T_2+R_1R_2T_3+U_1U_2T_4.
$$
Then
$$
T_4(n)C_{n+4}+T_3(n)C_{n+3}+T_2(n)C_{n+2}
-W_1(n)C_{n+1}-W_0(n)C_n=0.
$$

When $A_n=B_n$, the method has to be modified. Assume
$$
A_{n+2}=P(n)A_{n+1}+Q(n)A_n,
$$
where $P(n)$ and $Q(n)$ are rational. Then
$$
A_{n+3}=R(n)A_{n+1}+S(n)A_n,
$$
where
$$
R(n)=P(n+1)P(n)+Q(n+1), \quad S(n)=P(n+1)Q(n).
$$
Letting $C_n=A_n^2$ and $T(n)=R(n)Q(n)-P(n)S(n)$, we see that
$$
P(n)Q(n)C_{n+3}-R(n)S(n)C_{n+2}
-P(n)R(n)T(n)C_{n+1}+Q(n)S(n)T(n)C_n=0.
$$

The procedure leads usually (in cases (a)--(h)) to
a 8th order differential equation (6th order if $A_n=B_n$)
that factors into a 4th order equation. Cases \#167--179
were computed by this way.
\endremark

\medskip
Actually we found four 2nd order differential equations
somewhat suitable for Hadamard products:
$$
\allowdisplaybreaks
\xxalignat2
&\text{(k)} &
A_n&=9^n\sum_k(-1)^k\binom{-1/3}k\binom{-2/3}{n-k}\binom nk, \quad
D=\theta^2-3z(2\theta+1)-81z^2(\theta+1)^2;
\\
&\text{(l)} &
A_n&=8^n\sum_k(-1)^k\binom{-1/4}k\binom{-3/4}{n-k}\binom nk, \quad
D=\theta^2-4z(2\theta+1)-64z^2(\theta+1)^2;
\\
&\text{(m)} &
A_n&=36^n\sum_k(-1)^k\binom{-1/6}k\binom{-5/6}{n-k}\binom nk,
\\ &&
D&=\theta^2-24z(2\theta+1)-1296z^2(\theta+1)^2;
\\
&\text{(n)} &
A_n&=4^n\sum_k(-1)^k\binom{-1/2}k\binom{-1/2}{n-k}\binom nk, \quad
D=\theta^2-16z^2(\theta+1)^2.
\endxxalignat
$$
Only (n) gives Calabi-Yau differential equations when we
form $\text{(a)}*\text{(n)}$, $\text{(b)}*\text{(n)}$, $\dots$,
$\text{(n)}*\text{(n)}$ but this is just the substitution
$z\mapsto z^{1/2}$.
The Hadamard squares of (k), (l), (m) correspond to
Calabi--Yau differential equations found from the formula
$$
\align
&
(\theta^2-Az(2\theta+1)-Bz^2(\theta+1)^2)^{*2}
\\ &\quad
=\theta^4-z(2B\theta^4+A^2(2\theta+1)^2)
-Bz^2(2\theta+1)(2B\theta^2+2B\theta+4A^2+B)
\\ &\quad\qquad
+B^2z^3(2B\theta^4+8B\theta^3+(4A^2+12B)\theta^2+(4A^2+B)\theta+A^2+2B)
\\ &\quad\qquad
-B^4z^4(\theta+1)^4.
\endalign
$$

As a curiosity we mention that in the case $\text{(m)}*\text{(m)}$
we have
$$
q=z(1+21z+3840z^2+\dotsb)^{12960},
$$
where the series in the brackets is expected to lie in $\Bbb Z[[z]]$.

Of all non-square Hadamard products only $\text{(k)}*\text{(m)}$ gives
a Calabi--Yau equation
$$
\align
&
\theta^4-4\cdot18z(2\theta+1)^2
-2\cdot18^3z^2(27\theta^4+36\theta^3+74\theta^2+76\theta+24)
\\ &\quad
-32\cdot18^5z^3(2\theta+1)
+2\cdot18^7z^4(27\theta^4+72\theta^3+128\theta^2+72\theta+13)
\\ &\quad
+4\cdot18^9z^5(2\theta+1)^2-18^{12}z^6(\theta+1)^4.
\endalign
$$
To get integer instanton numbers we need to make the
substitution
$$
\wt K(q)
=K(q/3)=1+\sum_{k=1}^\infty\frac{k^3\wt N_kq^k}{1-q^k};
$$
then
$$
\wt N_1=-48, \quad
\wt N_2=-126, \quad
\wt N_3=-2864, \quad
\wt N_4=77958, \quad
\wt N_5=4942032, \quad \dotsc.
$$

\medskip
We can also get $5$th order equations by taking the Hadamard
product of (A)--(D) with $3$rd order equations of the following type:
$$
\allowdisplaybreaks
\xxalignat2
&(\alpha) &
A_n&=\sum_{k=0}^n\binom nk^2\binom{2k}k\binom{2n-2k}{n-k},
\\ &&
D&=\theta^3-2z(2\theta+1)(5\theta^2+5\theta+2)+64z^2(\theta+1)^3;
\\
&(\beta) &
A_n&=\sum_{k=0}^n\binom{2k}k^2\binom{2n-2k}{n-k}^2,
\\ &&
D&=\theta^3-8z(2\theta+1)(2\theta^2+2\theta+1)+256z^2(\theta+1)^3;
\\
&(\gamma) &
A_n&=\sum_{k=0}^n\binom nk^2\binom{n+k}k^2,
\\ &&
D&=\theta^3-z(2\theta+1)(17\theta^2+17\theta+5)+z^2(\theta+1)^3;
\\
&(\delta) &
A_n&=\sum_{k=0}^{\[n/3\]}(-1)^k3^{n-3k}\binom n{3k}
\binom{n+k}k\frac{(3k)!}{k!^3},
\\ &&
D&=\theta^3-z(2\theta+1)(7\theta^2+7\theta+3)+81z^2(\theta+1)^3;
\\
&(\epsilon) &
A_n&=\sum_k\binom nk^2\binom{2k}n^2,
\\ &&
D&=\theta^3-4z(2\theta+1)(3\theta^2+3\theta+1)+16z^2(\theta+1)^3;
\\
&(\zeta) &
D&=\theta^3-3z(2\theta+1)(3\theta^2+3\theta+1)-27z^2(\theta+1)^3;
\\
&(\eta) &
D&=\theta^3-z(2\theta+1)(11\theta^2+11\theta+5)+125z^2(\theta+1)^3,
\\
&(\vartheta) &
A_n&=(-64)^n\sum_k\binom{-1/2}k\binom{-1/2}{n-k}^3
=64^n\sum_k\binom{-1/4}k^2\binom{-3/4}{n-k}^2,
\kern20mm
\\ &&
D&=\theta^3-8z(2\theta+1)(8\theta^2+8\theta+5)+4096z^2(\theta+1)^3,
\\
&(\iota) &
A_n&=27^n\sum_k\binom{-1/3}k^2\binom{-2/3}{n-k}^2,
\\ &&
D&=\theta^3-3z(2\theta+1)(9\theta^2+9\theta+5)+729z^2(\theta+1)^3,
\\
&(\kappa) &
A_n&=432^n\sum_k\binom{-1/6}k^2\binom{-5/6}{n-k}^2,
\\ &&
D&=\theta^3-24z(2\theta+1)(18\theta^2+18\theta+13)+186624z^2(\theta+1)^3.
\endxxalignat
$$
The last 6 examples were found in~\cite{ES} (without formulas
for~$A_n$).

The general formula for the Hadamard product in this case is
as follows:
$$
\align
&
\bigl(\theta^2-z(2\theta+1)P(\theta)-cz^2(\theta+1)^3\bigr)
*\bigl(\theta^2-zQ(\theta)\bigr)
\\ &\qquad
=\theta^5-z(2\theta+1)P(\theta)Q(\theta)
-cz^2(\theta+1)Q(\theta+1)Q(\theta)
\endalign
$$
and we obtain the following table of corresponding cases%
\footnote{We computed all the other pullbacks of the Hadamard
products (one of those has trivial $K(q)$) and found no
new instanton numbers: they are all identical with earlier
known Hadamard products of type~$(a)*(b)$ etc.
The pullbacks are contained in the database~\cite{En}.}:
\medskip
\line{\hss\vbox{\offinterlineskip
\halign to55mm{\strut\tabskip=100pt minus 100pt
\strut\vrule\vphantom{\vrule height10.1pt}#&%
\hbox to10.5mm{\hfil$#$\hfil}&%
\vrule#&\hbox to10.5mm{\hfil#\hfil}&%
\vrule#&\hbox to10.5mm{\hfil#\hfil}&%
\vrule#&\hbox to10.5mm{\hfil#\hfil}&%
\vrule#&\hbox to10.5mm{\hfil#\hfil}&%
\vrule#\tabskip=0pt\cr\noalign{\hrule}
&          && (A) && (B) && (C) && (D) &\cr
\noalign{\hrule\vskip1pt\hrule}
& (\alpha) &&  39 &&  61 &&  37 &&  66 &\cr
& (\beta)  &&  40 &&  49 &&  43 &&  67 &\cr
& (\gamma) &&  44 &&  53 &&  52 && 149 &\cr
& (\delta) && 150 && 151 && 152 && 153 &\cr
\noalign{\hrule}
}}\hss}
\medskip

In addition, we have another type of Hadamard products.
If $\{A_n\}$ is given by any of the cases
($\alpha$)--($\epsilon$) with equation
$$
\theta^3-z(2\theta+1)P(\theta)+cz^2(\theta+1)^3,
$$
then the generating series for the sequence $\{\binom{2n}nA_n\}$
is annihilated by the 4th-order differential equation
$$
\theta^4-2z(2\theta+1)^2P(\theta)
+4cz^2(\theta+1)^2(2\theta+1)(2\theta+3).
$$
Thus we get cases \#16, 35, 29, 41, and 42, respectively.

\proclaim{Proposition 9}
Any differential equation of type $D'y=0$, where
$$
D'=\theta^3-z(2\theta+1)(a\theta^2+a\theta+b)+cz^2(\theta+1)^3,
$$
is the symmetric square of $Dy=0$, where
$$
D=\theta^2-z(2a\theta^2+a\theta+\tfrac b2)+cz^2(\theta+\tfrac12)^2.
\tag7.2
$$
\endproclaim

\demo{Proof}
In \cite{Al2} it is proved that the differential equation
$$
y^{(3)}+s_2(z)y''+s_1(z)y'+s_0(z)y=0
\tag7.3
$$
is a symmetric square of
$$
u''+p_1(z)u'+p_0(z)u=0
\tag7.4
$$
if and only if
$$
\frac{s_1s_2}3-\frac{2s_2^3}{27}+\frac{s_1'}2
-\frac{s_2''}6-\frac{s_2s_2'}3-s_0=0,
$$
and then
$$
p_1=\frac{s_2}3,
\qquad
p_0=\frac{s_1}4-\frac{s_2^2}{18}-\frac{s_2'}{12}.
$$
Here we get
$$
s_0=\frac{-b+cz}{1-2az+cz^2},
\qquad
s_1=\frac{1-(6a+2b)z+7cz^2}{1-2az+cz^2},
\qquad
s_2=\frac{3(1-3az+2cz^2)}{1-2az+cz^2},
$$
which satisfy the identity above. We also get
$$
p_1=\frac{1-3az+2cz^2}{z(1-2az+cz^2)},
\qquad
p_0=\frac{-2b+cz}{4z(1-2az+cz^2)},
$$
and the operator $\d^2/\d z^2+p_1(z)\,\d/\d z+p_0(z)$
is easily converted to~\thetag{7.2}.
\qed
\enddemo

\remark{Remark}
The Frobenius basis of the 3rd order differential
equation~\thetag{7.3} is
$$
y_0=u_0^2, \quad y_1=u_0u_1, \quad y_2=\frac12u_1^2,
$$
where $u_0,u_1$ form the Frobenius basis of~\thetag{7.4}.
Hence
$$
\frac{y_2}{y_0}=\frac12\biggl(\frac{u_1}{u_0}\biggr)^2
=\frac12\biggl(\frac{y_1}{y_0}\biggr)^2=\frac{t^2}2
$$
and
$$
K(q)=\frac{\d^2}{\d t^2}\biggl(\frac{y_2}{y_0}\biggr)=1
=\text{constant}.
$$

The substitution $z\mapsto 4z$ makes also
$u_0=y_0^{1/2}$ to have integer coefficients.
Since $\theta$ is invariant under this substitution we get
$$
D\mapsto\theta^2-z(8a\theta^2+4a\theta+2b)+4cz^2(2\theta+1)^2.
$$
Finding such equations is equivalent to finding integer solutions
to the recursion
$$
(n+1)^2u_{n+1}=2(a\cdot2n(2n+1)+2b)u_n-4c(2n-1)^2u_{n-1}
$$
which is similar to Zagier's computation in~\cite{Za}.

On the other hand, we found that only 2nd-order differential
operators of type
$\theta^2-zP(\theta)+cz^2(\theta+1)^2$,
where $P(\theta)$ is a polynomial of degree~2,
were useful for producing 4th-order Calabi--Yau
by taking the Hadamard product of two of them.
This, however, coincides
with Zagier's list~\cite{Za} of integer solutions to
$$
(n+1)^2A_{n+1}-an(n+1)A_n+bn^2A_{n-1}=\lambda A_n.
$$

Zagier's list has 36 entries. It contains our cases
(a)--(h), (A)--(C) (but not (D) since the coefficients
are too large). If $A_n$ is a polynomial in~$n$
(that corresponds to 8 entries in the list),
then the corresponding 2nd-order differential equation factors.
Excluding these polynomial cases, hypergeometric cases (A)--(C)
and terminating cases, we found that certain remaining
cases contain examples that we do not really like.
For instance, the recurrence
$$
(n+1)^2A_{n+1}-32n(n+1)A_n+256A_{n-1}=28A_n
$$
admits the following integer-valued solution:
$$
A_n=2^{4n}\sum_{k=0}^n\binom{1/2}{n-k}^2\binom{k+1/2}k
=\sum_{k=0}^n2^{2k}\frac{2k+1}{(2n-2k+1)^2}\binom{2n-2k}{n-k}^2\binom{2k}k;
$$
the 4th-order differential operator annihilating
the series $\sum_{n=0}^\infty\binom{2n}n^2A_nz^n$ is as follows:
$$
\theta^4-16z(2\theta+1)^2(8\theta^2+8\theta+7)
+4096z^2(2\theta+1)^2(2\theta+3)^2
$$
and the corresponding differential equation satisfies
relation~\thetag{2.2}; then computations show that
$z(q)\notin\Bbb Z[[q]]$ (and $N_l$ in~\thetag{2.9}
are not integers). This is an example
promised in Remark of Section~5.
\endremark

\head
8. More transformations
\endhead

{\it Transformations of the first $14$ hypergeometric cases}.
The success with negative rational numbers in binomial
coefficients caused us to make the following experiment.
Taking the hypergeometric operator in case~\#9,
$$
\theta^4-12^6z\biggl(\theta+\frac1{12}\biggr)\biggl(\theta+\frac5{12}\biggr)
\biggl(\theta+\frac7{12}\biggr)\biggl(\theta+\frac{11}{12}\biggr),
$$
consider
$$
A_n=1728^n\binom{2n}n\sum_k\binom{-1/12}k\binom{-5/12}k
\binom{-7/12}{n-k}\binom{-11/12}{n-k},
\qquad n=0,1,2,\dots.
$$
Then the series
$y_0(z)=\sum_{n=0}^\infty A_nz^n$
satisfies the differential equation $Dy=0$, where
$$
D=\theta^4-48z(2\theta+1)^2(72\theta^2+72\theta+41)
+2^{14}\cdot3^4z^2(2\theta+1)(2\theta+3)(3\theta+2)(3\theta+4),
$$
an equation found in~\cite{ES} and mentioned as case \#$9^*$
in \cite{AESZ, Table~A}. We get the instanton numbers
$$
\gathered
N_1=-480, \quad
N_2=-226848, \quad
N_3=16034720, \\
N_4=1094330202744, \quad
N_5=4352645747040, \quad \dots.
\endgathered
$$
We also note that
$$
y_0=(1+246z+768132z^2+\dotsb)^8,
\qquad
q=z(1+12z+33333z^2+\dotsb)^{288},
$$
where the series in brackets are expectively in $\Bbb Z[[z]]$.
Interchanging $5/12$ and $7/12$, we obtain case~$9^{**}$:
$$
\align
A_n
&=1728^n\binom{2n}n\sum_k\binom{-1/12}k\binom{-7/12}k
\binom{-5/12}{n-k}\binom{-11/12}{n-k}
\\
&=432^n\binom{2n}n^2\sum_k(-1)^k\binom{-5/6}k\binom{-1/6}{n-k}^2,
\endalign
$$
with the corresponding differential operator
$$
D=\theta^4-48z(2\theta+1)^2(72\theta^2+72\theta+31)
+2^{12}\cdot3^6z^2(2\theta+1)^2(2\theta+3)^2;
$$
the coupling $K(q)$ is the same as in \#$9^*$.

Trying the above experiment for all equations \#1--\#14
we find that cases \#1, \#2, \#11, \#12 and \#14
give integer mirror map and Yukawa coupling
but {\it there seems to be no way to get integer
instanton numbers\/} (by doing the substitution $z\to cz$).
Case \#$5^*$ has trivial $K(q)$. We have, however, a couple
of `good' cases \#$3^*$, \#$4^*$, \#$4^{**}$, \#$6^*$,
\#$7^*$, \#$7^{**}$, \#$8^*$, \#$8^{**}$, \#$10^*$, \#$10^{**}$,
\#$13^*$ and \#$13^{**}$ (see Appendix~A).
This (together with cases \#$9^*$, \#$9^{**}$ given above)
ends the list of the 14~ new equations.

By some unusual substitions we can also consider the following
additional cases.
In case \#$14^*$, we take
$$
A_n=144^n\binom{2n}n\sum_k\binom{-1/6}k\binom{-1/2}k
\binom{-5/6}{n-k}\binom{-1/2}{n-k}
$$
and obtain the corresponding differential operator
$$
D=\theta^4-2^4\cdot3^2z(2\theta+1)^2(2\theta^2+2\theta+1)
+2^{10}\cdot3^2z^2(2\theta+1)(2\theta+3)(3\theta+2)(3\theta+4).
$$
Then, choosing
$$
\wt K(q)=K\bigl((q/3)^{1/2}\bigr)
=1+\sum_{k=1}^\infty\frac{k^3\wt N_kq^k}{1-q^k},
$$
we have
$$
\gathered
\wt N_1=-2592, \quad
\wt N_2=-307800, \quad
\wt N_3=81451104, \\
\wt N_4=144135316512, \quad
\wt N_5=98667659422368, \quad \dotsc.
\endgathered
$$
In case \#$2^*$, when
$$
\gather
A_n=2000^n\binom{2n}n\sum_k\binom{-1/10}k\binom{-3/10}k
\binom{-7/10}{n-k}\binom{-9/10}{n-k},
\\
\aligned
D&
=\theta^4-2^4\cdot5z(2\theta+1)^2(50\theta^2+50\theta+33)
\\ \vspace{-2pt} &\qquad
+2^{10}\cdot5^4z^2(2\theta+1)(2\theta+3)(5\theta+4)(5\theta+6),
\endaligned
\endgather
$$
we need even more drastic methods;
the 'instanton numbers' land in $\Bbb Z[\sqrt5]$.
Taking
$$
\wt K(q)=K(q/\sqrt5)
=1+\sum_{k=1}^\infty\frac{k^3\wt N_kq^k}{1-q^k},
$$
we obtain
$$
\gathered
\wt N_1=-256\sqrt5, \quad
\wt N_2=32\sqrt5-35260, \quad
\wt N_3=1004288\sqrt5, \\
\wt N_4=835297220, \quad
\wt N_5=102454248704\sqrt5, \quad \dotsc.
\endgathered
$$
Case \#$2^{**}$ corresponds to the interchange
of $3/10$ and $7/10$, hence the change of the factor
$50\theta^2+50\theta+33$ by $50\theta^2+50\theta+17$
in the differential operator and conjugate $\wt K(q)$.
Is there a geometric interpretation of the numbers~$\wt N_k$?

\medskip
{\it Differential equations inspired by Guillera's formulas}.
The formulas discovered in~\cite{Gu}
caused us to try the following (again we mimic case \#9).
Let
$$
A_n'=3456^{2n}\sum_{k=0}^n(-1)^k\binom nk
\frac{(1/2)_k(1/12)_k(5/12)_k(7/12)_k(11/12)_k}{k!^5}.
$$
Then
$w_0=\sum_{n=0}^\infty A_n'z^n$
satisfies $D'w=0$, where
$$
\align
D'
&=\theta^5
-288z(2\theta+1)(103680\theta^4+207360\theta^3+262944\theta^2+159264\theta+41087)
\\ &\quad
+2^{20}3^8z^2(\theta+1)(207360\theta^4+829440\theta^3+1514592\theta^2+1370304\theta+498143)
\\ &\quad
-2^{38}3^{17}z^3(\theta+1)(\theta+2)(2\theta+3)(240\theta^2+720\theta+793)
\\ &\quad
+2^{53}3^{22}z^4(\theta+1)(\theta+2)(\theta+3)(360\theta^2+1440\theta+1633)
\\ &\quad
-2^{69}3^{30}z^5(\theta+1)(\theta+2)(\theta+3)(\theta+4)(2\theta+5).
\endalign
$$
The 4th order pullback of the latter differential operator is
$$
\align
D
&=\theta^4-2^43^2z(248832\theta^4+414720\theta^3+318528\theta
^2+111168\theta+14497)
\\ &\quad
+2^{22}3^{10}z^2(4\theta+3)(432\theta^3+1116\theta^2+886\theta+165)
\\ &\quad
-2^{34}3^{18}z^3(4\theta+1)(4\theta+3)(4\theta+7)(4\theta+9)
\endalign
$$
and, remarkably, it has lower degree than the 5th order
equation (cf\. case~\#129). We get the instanton numbers
$$
\gather
N_1=2710944,
\;\;
N_2=-717640978896,
\;\;
N_3=302270555492914464,
\\
N_4=-171507700573958028578832,
\;\;
N_5=113303073680022744870130144224,
\;\; \dots,
\endgather
$$
and note that
$$
\gather
y_0=(1+260946z+1405445560884z^2+\dotsb)^8,
\\
q=x(1+13295z+67507583411z^2+\dotsb)^{576},
\quad
K(q)=(1+56478q-\dotsb)^{48},
\endgather
$$
where the parentheses are supposed to have integer coefficients.

Doing the same thing with case \#3, i.e\. taking
$$
A_n'=1024^n\sum_{k=0}^n(-1)^k\binom nk\frac{(1/2)_k^5}{k!^5},
$$
we get the same Yukawa coupling as for case \#115 which is an Hadamard
square with
$$
A_n=\biggl\{\sum_k4^{n-k}\binom{2k}k^2\binom{2n-2k}{n-k}^2\biggr\}^2
$$
But the differential equations are quite different as are
the solutions~$y_0$.

We list all possible variations on the theme
as cases \#$\wh1$--$\wh{14}$ in \cite{AESZ, Table~A}.

\smallskip
There is an even simpler way to arrive at essentially the same
Yukawa coupling $K(q)$ as above.
Consider the 5th order hypergeometric operator
$$
\theta^5-4\cdot12^6z\Bigl(\theta+\frac12\Bigr)\Bigl(\theta+\frac1{12}\Bigr)
\Bigl(\theta+\frac5{12}\Bigr)\Bigl(\theta+\frac7{12}\Bigr)
\Bigl(\theta+\frac{11}{12}\Bigr),
$$
with the solution
$$
w_0=\sum_{n=0}^\infty
\frac{(12n)!}{(6n)!n!^6}\frac{z^n}{\binom{4n}{2n}}
$$
of the corresponding differential equation $D'w=0$, and
with the 4th order pullback
$$
\align
D
&=\theta^4
-2^43^2z(331776\theta^4+82944\theta^3+13248\theta^2-28224\theta-14497)
\\ &\quad
+2^{19}3^8z^2(248832\theta^4+124416\theta^3+25056\theta^2-4176\theta+21143)
\\ &\quad
-2^{32}3^{14}z^3(331776\theta^4+248832\theta^3+60480\theta^2+21312\theta-4453)
\\ &\quad
+2^{51}3^{22}z^4\theta(2\theta+1)(12\theta+1)(12\theta+5)
\endalign
$$
Observe that here the pullback increased the degree from~1 to~4.
The instanton numbers are
$$
N_1=-2710944,
\quad
N_2=-717640301160,
\quad \dots,
$$
i.e\. the same as above up to a sign. We also note
$$
y_0=(1-86982z-577991455848z^2-\dotsb)^{24},
\;\;
q=x(1+7441z+32486988467z^2+\dotsb)^{576}.
$$

A third variation of the theme is the following. Consider
$$
A_n'=(4\cdot 12^6)^n\binom{-1/2}n\sum_{k=0}^n\binom nk
\binom{-1/12}k\binom{-5/12}k\binom{-7/12}k\binom{-11/12}k
$$
with corresponding 5th order differential operator
$$
\align
D'
&=\theta^5
+2\cdot 12^2z(2\theta+1)(124416\theta^4+248832\theta^3+234720\theta^2+110304\theta+21121)
\\ &\quad
+112\cdot 12^{10}z^2(\theta+1)(2\theta+1)(2\theta+3)(72\theta^2+144\theta+101)
\\ &\quad
+16\cdot 12^{16}z^3(2\theta+1)(2\theta+3)(2+5)(1152\theta^2+3456\theta+2831)
\\ &\quad
+12^{26}z^4(\theta+2)(2\theta+1)(2\theta+3)(2\theta+5)(2\theta+7)
\\ &\quad
+64\cdot 12^{30}z^5(2\theta+1)(2\theta+3)(2\theta+5)(2+7)(2\theta+9)
\endalign
$$
with pullback
$$
\align
D
&=\theta^4
+2^43^2z(995328\theta^4+497664\theta^3+220608\theta^2-28224\theta-35233)
\\ &\quad
+2^{18}3^8z^2(5142528\theta^4+5308416\theta^3+2946816\theta^2+9792\theta-292949)
\\ &\quad
+2^{35}3^{14}z^3(1866240\theta^4+2985984\theta^3+2011320\theta^2+142920\theta-208501)
\\ &\quad
+2^{48}3^{20}z^4(6656256\theta^4+14681088\theta^3+11732832\theta^2+1326960\theta-1309273)
\\ &\quad
+2^{64}3^{28}z^5(4\theta-1)(51840\theta^3+160704\theta^2+177908\theta+64537)
\\ &\quad
+2^{79}3^{34}z^6(4\theta-1)(4\theta+3)(4464\theta^2+13536\theta+10985)
\\ &\quad
+2^{97}3^{43}z^7(\theta+2)(4\theta-1)(4\theta+3)(4\theta+7)
\\ &\quad
+2^{108}3^{48}z^8(4\theta-1)(4\theta+3)(4\theta+7)(4\theta+11).
\endalign
$$
We have
$$
y_0=(1+211398z-2341345644648z^2+\dotsb)^{24}
$$
and
$$
q=z(1-28177z+245565832115z^2-\dotsb)^{576}.
$$
The instanton numbers are the same as above. Thus we have three fourth order
differential equations of degree 3, 4, and 8, respectively.
The instanton numbers, however, are invariants
(possibly, with some geometric interpretation). It
would be interesting to find the transformations
between the three different solutions.

\medskip
{\it The mirror at infinity}.
In \cite{Ro} E.~R\o dland studied case \#27 at infinity.
We choose instead \#124 from~\cite{ES}:
$$
\align
D
&=61^2\theta^4
-61z(3029\theta^4+5572\theta^3+4677\theta^2+1891\theta+305)
\\ &\qquad
+z^2(1215215\theta^4+3428132\theta^3+4267228\theta^2+2572675\theta+611586)
\\ &\qquad
-3^4z^3(39370\theta^4+140178\theta^3+206807\theta^2+142191\theta+37332)
\\ &\qquad
+3^8z^4(566\theta^4+2230\theta^3+3356\theta^2+2241\theta+558)
-3^{13}z^5(\theta+1)^4,
\endalign
$$
where the explicit formulas for $A_n$ are known
and one should take $N_0=61$ in order to get
integer instanton numbers (cf.\ \thetag{2.9}).
The substitution $z\mapsto3^{-5}z^{-1}$, $y\mapsto zy$
results in $\theta\to-\theta-1$
(this works only if the highest degree term has the form
$(\theta+1)^4$) and we obtain the `dual' differential operator
$$
\align
D^*
&=\theta^4-z(566\theta^4+34\theta^3+62\theta^2+45\theta+9)
\\ &\qquad
+3z^2(39370\theta^4+17302\theta^3+22493\theta^2+8369\theta+1140)
\\ &\qquad
-3^2z^3(1215215\theta^4+1432728\theta^3+1274122\theta^2+538245\theta+93222)
\\ &\qquad
+3^761z^4(3029\theta^4+6544\theta^3+6135\theta^2+2863\theta^2+548)
-3^{12}61^2z^5(\theta+1)^4.
\endalign
$$
We do not know a formula for the corresponding $A_n$ here.

Many differential equations with highest degree term
$(\theta+1)^4$ (in particular, all Hadamard products) are
self-dual. There is, however, the class with
$$
A_n=\sum_k\binom nk^{5-2r}\binom{2k}k^r\binom{2n-2k}{n-k}^r
$$
for $r=0,1,2,3,4,5$ (cases \#22, \#21, \#23, \#56, \#71, \#118,
respectively). Then the duality maps
$r\leftrightarrow5-r$.

Also differential equations with highest degree term
$(2\theta+1)^4$ can be reflected in~$\infty$. Let us take, e.g.,
case~\#55 with
$$
\align
D
&=9\theta^4-12z(208\theta^4+224\theta^3+163\theta^2+51\theta+6)
\\ &\qquad
+2^9z^2(32\theta^4-928\theta^3-1606\theta^2-837\theta-141)
\\ &\qquad
+2^{16}z^3(144\theta^2+576\theta^3+467\theta^2+144\theta+15)
-2^{24}z^4(2\theta+1)^4.
\endalign
$$
Then the transformation $z\to2^{-18}z^{-1}$ maps
$\theta$ to $-\theta-1/2$ and gives the dual
$$
\align
D^*
&=\theta^4-2^4z(576\theta^4-1152\theta^3-724\theta^2-148\theta-13)
\\ &\qquad
-2^{17}z^2(32\theta^4+992\theta^3-166\theta^2-57\theta-6)
\\ &\qquad
+2^{26}3z^3(832\theta^4+768\theta^3+556\theta^2+192\theta+25)
-2^{40}3^2z^4(2\theta+1)^4
\endalign
$$
Also cases \#33, \#99, \#154 can be treated in the same way
(case \#154 is self-dual).

\head
9. Supercongruences and $k$-realizable series
\endhead

It is time to stop and discuss the phenomenon of
integer coefficients of Yukawa couplings number-theoretically.

Let $T\:X\to X$ be a map of a set. Define the number
of points~$x$ with period~$n$:
$$
f_n(T)=\#\{x\in X:T^nx=x\},
$$
and the number of points $x$ with the smallest period~$n$:
$$
f_n^*(T)=\#\bigl\{x\in X:T^nx=x,\ \#\{T^lx:l\in\Bbb Z\}=n\bigr\}.
$$
Then $f_n^0(T)=\frac1nf_n^*(T)$ determines the number of orbits
of length~$n$. Applying the M\"obius inversion formula
to $f_n(T)=\sum_{d\mid n}f_d^*(T)$ we obtain
$$
f_n^0(T)=\frac1n\sum_{d\mid n}\mu\Bigl(\frac nd\Bigr)f_d(T).
\tag9.1
$$

\proclaim{Proposition 10 \cite{PW}}
Let $\{A_n\}_{n=1}^\infty$ be a sequence with all $a_n\ge0$.
Then, there exists a map $T\:X\to X$ such that $A_n=f_n(T)$
if and only if the numbers
$$
B_n=\frac1n\sum_{d\mid n}\mu\Bigl(\frac nd\Bigr)A_d,
\qquad n=1,2,\dots,
$$
are all non-negative integers.
\endproclaim

\demo{Proof}
The `$\!\implies\!$'-part was already done in~\thetag{9.1}.

In order to prove the `$\!\impliedby\!$'-part, we take
$X=\{0,1,2,\dots\}$ and let $T\:X\to X$~be a permutation
consisting of~$B_n$ cycles of length~$n$ for all $n=1,2,\dots$\,.
Then $A_n=f_n(T)$ as required.
\qed
\enddemo

To illustrate the above construction of~$T$ in the proof,
we consider the following

\example{Example}
Let $\{A_n\}_{n=1}^\infty=\{1,5,7,17,31,65,\dots\}$ be given
by the recurrence
$A_{n+1}=\allowmathbreak A_n+2A_{n-1}$ for $n=1,2,\dots$\,.
Then $\{B_n\}_{n=1}^\infty=\{1,2,2,3,6,9,\dots\}$ and a choice
of~$T$ is as follows:
$$
T=\underbrace{(0)}_{B_1=1}\,
\underbrace{(1 \; 2)(3 \; 4)}_{B_2=2}\,
\underbrace{(5 \; 6 \; 7)(8 \; 9 \; 10)}_{B_3=2}\,
\underbrace{(11 \; 12 \; 13 \; 14)(15 \; 16 \; 17 \; 18)
(19 \; 20 \; 21 \; 22)}_{B_4=3}\dots\,.
$$
The points of order $4$ are
$$
\{11,12,\dots,22, \; 1,2,3,4, \; 0\},
$$
hence $f_4(T)=4\cdot3+2\cdot2+1=17=A_4$.
\endexample

Assume now that for some $k\ge1$ we have
$$
\frac1{n^k}\sum_{d\mid n}\mu\Bigl(\frac nd\Bigr)A_d\in\Bbb N
\qquad\text{for all $n=1,2,\dots$}\,.
$$
This means that the number of orbits of length $n$
is divisible by~$n^{k-1}$.
The sequence $\{A_n\}_{n=1,2,\dots}$ is said to be
$k$-{\it realizable\/} if
$$
B_n=\frac1{n^k}\sum_{d\mid n}\mu\Bigl(\frac nd\Bigr)A_d\in\Bbb Z
\qquad\text{for all $n=1,2,\dots$}
\tag9.2
$$
holds. Call a power series $y=\sum_{n=0}^\infty A_nz^n\in\Bbb Q[[z]]$
\ $k$-{\it realizable\/} if there exists an integer $C$ such that
$\{CA_n\}_{n=1}^\infty$ is $k$-realizable.

The most remarkable property of the examples of the Calabi--Yau
differential equations is that the Yukawa coupling~\thetag{2.9}
is (expected to be) $3$-realizable; the pseudo-coupling~\thetag{3.8}
constructed in Section~3 seems to be $2$-realizable (where in cases
\#15--29 we have to replace $\Bbb Z$ in~\thetag{9.2} by $\Bbb Z[1/p]$
for a suitable choice of the prime~$p$). Another striking property
of the first 14 (hypergeometric) cases from \cite{AESZ, Table~A} is that
$y_0(z)$ is also $3$-realizable.

The property of being $k$-realizable is related to
Kummer supercongruences by the following result.

\proclaim{Proposition 11}
The following are equivalent:
\roster
\item"(i)" $\{A_n\}$ is $k$-realizable;
\item"(ii)" $A_{mp^r}\equiv A_{mp^{r-1}}\pmod{p^{rk}}$
for $p$~prime and $(m,p)=1$;
\item"(iii)"
$\sum_{n=1}^\infty A_nz^n/n^k=\sum_{n=1}^\infty B_n\Li_k(z^n)$,
where $B_n$ are the integers defined in~\thetag{9.2} and
$$
\Li_k(x)=\sum_{n=1}^\infty\frac{x^n}{n^k}
$$
denotes the $k$-th polylogarithm.
\endroster
\endproclaim

\demo{Proof}
The proofs on pp.~3--6 in \cite{Re} for $k=1$ are easily
generalized; note that $\Li_1(x)=-\log(1-x)$.
\qed
\enddemo

\remark{Remark}
For the first $n$ numbers $A_1,A_2,\dots,A_n$, the case $k=1$
is well described in Exercise 5.2 of~\cite{St}. Then there is
also an ($n\times n$)-matrix $M$ with integer entries such that
$A_k=\Tr(M^k)$ for $k=1,2,\dots,n$ (see~\cite{Al1}).
\endremark

\head
10. $6$th- and higher order differential equations
\endhead

We would like to start this last section by presenting two other
puzzling examples related to 4-term polynomial recursions
$$
\align
&
(n+1)^6(41218n^3-48459n^2+20010n-2871)A_{n+1}
\\ &\quad
+2(48802112n^9+89030880n^8+36002654n^7-24317344n^6-19538418n^5
\\ &\quad\qquad
+1311365n^4+3790503n^3+460056n^2-271701n-60291)A_n
\\ &\quad
-4(2n-1)
(3874492n^8-2617900n^7-3144314n^6+2947148n^5+647130n^4
\\ &\quad\qquad
-1182926n^3+115771n^2+170716n-44541)A_{n-1}
\\ &\quad
-4(n-1)^4(2n-1)(2n-3)(41218n^3+75195n^2+46746n+9898)A_{n-2}=0
\tag10.1
\endalign
$$
for simultaneous approximations to $\zeta(3)$ and $\zeta(5)$
(see \cite{Zu2}), and
$$
\align
&
n(n+1)^5(91n^3-182n^2+126n-30)A_{n+1}
\\ &\quad
-n(3458n^8+1729n^7-2947n^6-2295n^5+901n^4+1190n^3
\\ &\quad\qquad
+52n^2-228n-60)A_n
\\ &\quad
-(153881n^9-307762n^8+185311n^7+2960n^6-31631n^5-88n^4
\\ &\quad\qquad
+5239n^3-610n^2-440n+100)A_{n-1}
\\ &\quad
+24(n-1)^3(2n-1)(6n-5)(6n-7)(91n^3+91n^2+35n+5)A_{n-2}=0
\tag10.2
\endalign
$$
for the sequence $\sum_{k=0}^n\binom nk^6$
(see \cite{Pe}, \cite{SJ}).

The differential equations for generating series of \thetag{10.1}
and \thetag{10.2} are of order~9; nevertheless, {\tt Maple}
easily factorise the corresponding differential operators
and we get the 6th-order linear differential equations
$$
\align
&
z^5(196z-87)^3(16z^3+752z^2-2368z-1)\frac{\d^6}{\d z^6}
\\ &\quad
+3z^4(196z-87)^2(21952z^4+873216z^3-2778608z^2+1235312z+435)\frac{\d^5}{\d z^5}
\\ &\quad
+z^3(196z-87)(85898176z^5+2768881024z^4-8828169756z^3
\\ &\quad\qquad
+7144975624z^2-1768825884z-491985)\frac{\d^4}{\d z^4}
\allowdisplaybreak &\quad
+2z^2(20932110080z^6+508609400320z^5-1613572776144z^4+1739040695000z^3
\\ &\quad\qquad
-824264516904z^2+148195933632z+29632635)\frac{\d^3}{\d z^3}
\allowdisplaybreak &\quad
+z(35509291776z^6+579712191744z^5-1530351585392z^4+1500993519824z^3
\\ &\quad\qquad
-731658297456z^2+173252093886z+20413593)\frac{\d^2}{\d z^2}
\allowdisplaybreak &\quad
+(7138000128z^6+55844570880z^5-123313425872z^4+96429989856z^3
\\ &\quad\qquad
-16021623504z^2+17983065996z+658503)\frac{\d}{\d z}
\allowdisplaybreak &\quad
-6(15059072z^5-4148928z^4+9924264z^3+214891044z^2
\\ &\quad\qquad
+106071966z+4609521)
=0
\tag10.3
\endalign
$$
and
$$
\align
&
z^4(z-1)(64z-1)(27z+1)(75z^3+1420z^2+561z+9)\frac{\d^6}{\d z^6}
\\ &\quad
+z^3(2786400z^6+52926750z^5-22883417z^4-19427551z^3
\\ &\quad\qquad
-654306z^2+1308z+126)\frac{\d^5}{\d z^5}
\allowdisplaybreak &\quad
+5z^2(3822480z^6+77258112z^5-18886036z^4-23855750z^3
\\ &\quad\qquad
-758207z^2-1275z+99)\frac{\d^4}{\d z^4}
\allowdisplaybreak &\quad
+5z(9906840z^6+214908768z^5-10697346z^4-53272456z^3
\\ &\quad\qquad
-1605389z^2-7047z+117)\frac{\d^3}{\d z^3}
\allowdisplaybreak &\quad
+2(22370400z^6+528293250z^5+82763885z^4-96071120z^3
\\ &\quad\qquad
-2836582z^2-15165z+72)\frac{\d^2}{\d z^2}
\allowdisplaybreak &\quad
+6(1684500z^5+44517700z^4+17102475z^3
\\ &\quad\qquad
-4687335z^2-158967z-540)\frac{\d}{\d z}
\\ &\quad
+12(15750z^4+503175z^3+327205z^2-1845z-1044)
=0,
\tag10.4
\endalign
$$
respectively. (The corresponding recursions will be
7-term but we do not use this fact.)

The differential equation~\thetag{10.3} is MUM, so taking
the Frobenius basis $y_0,y_1,y_2,\allowmathbreak\dots,y_5$ we may compute
the inverse of $q(z)=\exp(y_1/y_0)$,
$$
\align
z(q)
&=q+230q^2-26827q^3+24147708q^4-23334210874q^5+26542920855790q^6
\\ &\qquad
-33500728089853156q^7+45492345805504886104q^8+O(q^9),
\endalign
$$
and pseudo-coupling
$$
\tilde K(q)
=\biggl(q\frac{\d}{\d q}\biggr)^2\biggl(\frac{y_2}{y_0}\biggr)
=1+\sum_{l=1}^\infty\frac{l^2\tilde N_lq^l}{1-q^l}
\tag10.5
$$
where
$$
\gathered
\tilde N_1=-320, \quad
\tilde N_2=118264, \quad
\tilde N_3=-84117876, \quad
\tilde N_4=80349364184,
\\
\tilde N_5=-90632838175404, \quad
\tilde N_6=113783008482427048,
\\
\tilde N_7=-153937885949108788148, \quad
\tilde N_8=220092295805975113694144, \quad
\dots
\endgathered
$$
are expected to be all integers.

The case of the (not MUM) differential equation~\thetag{10.4} is
more delicate: the exponents at the point $z=0$ are $0,0,0,0,0,1$,
so that we get two linearly independent analytic solutions at the point.
The choice of the first three components of the basis at $z=0$ are
as follows:
$$
\align
y_0
&=1+2z+66z^2+1460z^3+54850z^4+2031252z^5+86874564z^6
\\ &\qquad
+3848298792z^7+180295263810z^8+O(z^9)
\allowdisplaybreak
y_1
&=y_0\log z+
6z+201z^2+5114z^3+\frac{392153}2z^4+\frac{37695381}5z^5
\\ &\qquad
+\frac{1636192577}5z^6+\frac{515198795302}{35}z^7
+\frac{19496105427567}{28}z^8+O(z^9),
\allowdisplaybreak
y_2
&=-y_0\frac{\log^2z}2+y_1\log z
+\beta\biggl(z+\frac{45}4z^2+\frac{7085}{18}z^3+\frac{1750325}{144}z^4
+\frac{293631899}{600}z^5
\\ &\qquad
+\frac{4039916881}{200}z^6
+\frac{4451078944741}{4900}z^7+\frac{663063626523441}{15680}z^8\biggr)
\\ &\qquad
+108z^2+3690z^3+164379z^4+6851805z^5+\frac{3138305727}{10}z^6
\\ &\qquad
+\frac{73262810391}5z^7+\frac{99760681900977}{140}z^8+O(z^9)
\endalign
$$
(for the series $y_2$ we have an additional `freedom' $\beta\in\Bbb Q$).
Then we get the expansion
$$
z(q)=q-6q^2-135q^3-380q^4-24960q^5-696366q^6-26153302q^7
-901888104q^8+O(q^9)
$$
and the coefficients in the Lambert expansion~\thetag{10.5}
$$
\gathered
\tilde N_1=\beta, \quad
\tilde N_2=3(\beta+30), \quad
\tilde N_3=63(\beta+20), \quad
\tilde N_4=1357\beta+43200,
\\
\tilde N_5=5(9139\beta+235152), \quad
\tilde N_6=3(479743\beta+14248200),
\\
\tilde N_7=126(419883\beta+11883130), \quad
\tilde N_8=3(653810477\beta+19370535600), \quad
\dots
\endgathered
$$
are expected to be in $\Bbb Z\beta+\Bbb Z$.

We have several further examples of 6th-order differential
equations with the same collection of the exponents at the point
$z=0$ as for equation \thetag{10.4}. In all such examples we get
the similar phenomenon of admitting a free parameter~$\beta$
for the coefficients in the pseudo-coupling~\thetag{10.5}.
Here we only mention formulas for $\{A_n\}_{n=0,1,\dots}$
in the examples, where the series $y_0=\sum_{n=1}^\infty A_nz^n$
form analytic solutions in $\Bbb Z[[z]]$ to the corresponding
differential equations:
$$
A_n=\frac{(3n)!}{n!^3}\sum_{k=0}^n\binom nk^4,
\qquad
A_n=\frac{(4n)!}{n!^4}\sum_{j,k}\binom nj\binom nk^2\binom kj.
$$

Similar phenomena seem to be happened in more general situations.
For example, the series
$$
y_0(z)
=\sum_{n=0}^\infty z^n\frac{(4n)!}{n!^4}\sum_{k=0}^n\binom nk^4
$$
gives an analytic solution to the 8th-order linear differential
equation $Dy=0$, where
$$
\align
D
&=\theta^6(\theta-1)^2
-16z\theta^2(2\theta+1)^2(4\theta+1)(4\theta+3)(3\theta^2+3\theta+1)
\\ &\qquad
-256z^2(2\theta+1)(2\theta+3)
(4\theta+1)(4\theta+3)^2(4\theta+5)^2(4\theta+7);
\endalign
$$
a suitable choice of the solution $y_1(z)$ leads to
the inversion $z(q)$ of the series $q(z)=\allowmathbreak\exp(y_1/y_0)$
satisfying (conjecturally!) $z(q)\in\Bbb Z[[q]]$.
Finally, one can take a free
parameter $\beta$ for the solution $y_2(z)$
giving the pseudo-coupling~\thetag{10.5} with the coefficients
$$
\gathered
\tilde N_1=\beta, \quad
\tilde N_2=199\beta+308464, \quad
\tilde N_3=17(18671\beta+18181888),
\\
\tilde N_4=2(233056091\beta+290786122832), \
\tilde N_5=3(311572981529\beta+381333083404544),
\\
\tilde N_6=2020548109265033\beta+2585332849682835728, \quad
\dots
\endgathered
$$
that are expected to be in $\Bbb Z\beta+\Bbb Z$.
Another interesting example is the series
$$
y_0(z)
=\sum_{n=0}^\infty z^n\sum_{k=0}^n\binom nk^7,
$$
which also satisfies a (linearly irreducible)
8th-order differential equation $Dy=0$ with exponents
$0,0,0,0,0,0,1,1$ about $z=0$.


\vskip5mm
\bgroup
\eightpoint\smc

\hbox to70mm{\kern2mm\vbox{\hsize=70mm%
\leftline{Matematikcentrum}
\leftline{Lunds Universitet}
\leftline{Matematik MNF, Box 118}
\leftline{SE-221\,00 Lund, SWEDEN}
\leftline{{\it E-mail address\/}: \tt gert\@maths.lth.se}
}\hss}

\bigskip
\hbox to70mm{\kern2mm\vbox{\hsize=70mm%
\leftline{Department of Mechanics and Mathematics}
\leftline{Moscow Lomonosov State University}
\leftline{Vorobiovy Gory, GSP-2}
\leftline{119992 Moscow, RUSSIA}
\leftline{{\it E-mail address\/}: \tt wadim\@ips.ras.ru}
}\hss}

\egroup
\enddocument